\documentclass[mnsc]{informs1}              

\usepackage{natbib}
 \NatBibNumeric
 \def\bibfont{\small}%
 \bibpunct[, ]{[}{]}{,}{n}{}{,}%

\usepackage[linesnumbered,ruled]{algorithm2e}
\usepackage{tikz}
\usetikzlibrary{shapes,arrows}
\usepackage{pgfplots}
\usepackage{mathrsfs}
\usepackage[colorlinks=true,breaklinks=true,bookmarks=true,urlcolor=blue,
     citecolor=blue,linkcolor=blue,bookmarksopen=false,draft=false]{hyperref}

\def\EMAIL#1{\href{mailto:#1}{#1}}
\def\URL#1{\href{#1}{#1}}         

\TheoremsNumberedThrough     

\EquationsNumberedThrough    

\renewcommand*{\vec}[1]{\vbox{\halign{##\cr 
  \tiny\rightarrowfill\cr\noalign{\nointerlineskip\vskip1pt} 
  $#1\mskip2mu$\cr}}}

\newcommand{\eq}{\begin{equation}}
\newcommand{\qe}{\end{equation}}

\newcommand{\R}{\mathbb{R}}

\theoremstyle{plain}

\theoremstyle{definition}

\newcommand{\inte}{\displaystyle\int}

\newcommand{\integ}{\displaystyle\int_{\Gamma_G}}

\newcommand{\intev}{\displaystyle\int_{V}}
\newcommand{\intega}{\displaystyle\int_{\gamma_{x,y}}}

\newcommand{\intee}{\displaystyle\int_{e}}

\newcommand{\inteve}{\displaystyle\int_{V_e}}

\newcommand{\inteves}{\displaystyle\int_{V_{e_{\tilde{s}}}}}

\newcommand{\bule}{\mathbb{B}\left(x^\star,\dfrac{1}{\beta}\right)}
\newcommand{\de}{\mathrm{d}}

\newcommand{\dm}{\mathrm{d}m_t(x,y)}
\newcommand{\dn}{\mathrm{d}n_t(x)}
\newcommand{\dt}{\partial_t}
\newcommand{\dx}{\mathrm{d}x}
\newcommand{\dy}{\mathrm{d}y}
\newcommand{\dxy}{\mathrm{d}x\mathrm{d}y}

\newcommand{\mub}{\mu_{\beta_t}}
\newcommand{\ltd}{\tilde{\mathcal{L}}_{2,t}}
\newcommand{\lcd}{\hat{\mathcal{L}}_{2,t}}

\newcommand{\um}{U_{\nu}}
\newcommand{\rt}{\tilde{R}_t}
\newcommand{\rc}{\hat{R}_t}
\newcommand{\lu}{\mathcal{L}_{1,t}}
\newcommand{\ld}{\mathcal{L}_{2,t}}
\newcommand{\at}{\alpha_t}
\newcommand{\z}{Z_{\beta_t}}
\newcommand{\eb}{\mathcal{E}_{\beta}}

\newcommand{\zb}{Z_{\beta}}
\newcommand{\en}{\mathbf{Ent}_{\mu}}
\newcommand{\enb}{\mathbf{Ent}_{\mubb}}
\newcommand{\enbt}{\mathbf{Ent}_{\mub}}
\newcommand{\grad}{\triangledown}
\newcommand{\mubb}{\mu_{\beta}}
\newcommand{\lot}{\log \dfrac{n_t(x)}{\mub (x)}}
\newcommand{\fot}{f_t(x)}

\usepackage{color}
\usepackage{multicol}

\newcolumntype{M}[1]{>{\raggedright}m{#1}}
\begin{document}



\RUNTITLE{On the calculation of a graph barycenter}

\TITLE{How to calculate the barycenter of a weighted graph}

\ARTICLEAUTHORS{%
\AUTHOR{S\'ebastien Gadat}
\AFF{Toulouse School of Economics, Universit\'e Toulouse 1 Capitole \EMAIL{sebastien.gadat@math.univ-toulouse.fr}, \URL{http://perso.math.univ-toulouse.fr/gadat/}}
\AUTHOR{Ioana Gavra}
\AFF{Institut de Math\'ematiques de Toulouse, Universit\'e Toulouse 3 Paul Sabatier \EMAIL{Ioana.Gavra@math.univ-toulouse.fr}, \URL{}}
\AUTHOR{Laurent Risser}
\AFF{Institut de Math\'ematiques de Toulouse, CNRS
\EMAIL{lrisser@math.univ-toulouse.fr}, \URL{http://laurent.risser.free.fr/}}

} 

\ABSTRACT{%
Discrete structures like graphs make it possible to naturally and flexibly model complex phenomena. 
Since graphs that represent various types of information are increasingly available today, 
their analysis has become a popular subject of research. The graphs studied in the field of data science at this time generally have a 
large number of nodes that  are not fairly weighted and connected to each other, translating a structural specification of the data. Yet, even an algorithm for locating the \textit{average} position in graphs is lacking although this knowledge would be of primary interest for statistical or representation problems.
In this work, we develop a stochastic algorithm for finding the Fr\'echet mean of weighted undirected metric graphs. This method relies on a noisy simulated annealing algorithm dealt with using homogenization. We then illustrate our algorithm with two examples  (subgraphs  of a social network and of a collaboration and citation network).
}%

\MSCCLASS{Primary: 90B05; secondary: 90C40, 90C90}

\MSCCLASS{60J20,60J60,90B15}
\ORMSCLASS{Primary: Probability - Markov processes; secondary: Networks/graphs - Traveling salesman}
\HISTORY{}

\maketitle

%
\section{Introduction.}\label{sec:intro}


Numerous open questions in a very wide variety of scientific domains involve complex discrete structures that are easily modeled by graphs. The nature of these graphs may be weighted or not, directed or not, observed online or by using   batch processing, each time implying new problems and sometimes leading to difficult mathematical or numerical questions.
Graphs are the subject of perhaps one of the most impressive growing bodies of literature dealing with  potential applications in statistical or quantum physics (see, \textit{e.g.}, \cite{estrada}) economics (dynamics in economy structured as networks), biology (regulatory networks of genes, neural networks), informatics (Web understanding and representation), social sciences (dynamics in social networks, analysis of citation graphs).
We refer to \cite{kolaczyk} and \cite{newman} for recent communications on the theoretical aspects of random graph models, questions in the field of statistics and graphical models, and related numerical algorithms. 
In \cite{jackson},  the authors have developed an overview of numerous possible applications using graphs and networks in the fields of industrial organization and economics.  Additional applications, details and references in the field of machine learning may also be found in \cite{Goldenberg}.


In the meantime, the nature of the mathematical questions raised by the models that involve networks is  very extensive and may concern geometry, statistics, algorithms or dynamical evolution over the network, to name a few.
 For example, we can be interested in the definition of  suitable random graph models that make it possible to detect specific shape phenomena observed at different scales and frameworks 
 (the small world networks of \cite{Watts}, the existence of a giant connected component for a specific range of parameters in the Erd\"os-R\'enyi random graph model \cite{erdos}). A complete survey may be found in \cite{lovasz}.
 Another important field of investigation is dedicated to graph visualization (see some popular methods in   \cite{Shneiderman} and \cite{Kaufmann}, for example).
 In statistics, a popular topic deals with the estimation of a natural clustering when the graph follows a specific random graph model (see, among others, the recent contribution \cite{Tsybakov} that proposes an optimal estimator for the stochastic block model and  smooth graphons). Other approaches rely on efficient algorithms that analyze the spectral properties of adjacency matrices representing the networks (see, \textit{e.g}, \cite{bach}).
 A final important field of interest deals with the evolution of a dynamical system defined through a discrete graph structure: this is, for instance, the question raised by gossip models that may describe belief evolution over a social network (see \cite{acemoglu} and the references therein).


 We address a problem here that may be considered as very simple at  first glance: we aim to define and estimate the barycenter of weighted graphs. Hence, this problem involves questions that straddle the area of statistics and the geometry of graphs. Surprisingly, as far as we know, this question has received very little attention, although a good assessment of the location of a weighted graph barycenter could be used for fair representation issues or for  understanding the graph structure from a statistical point of view. In particular, it could be used and extended to produce ``second order'' moment analyses of graphs. We could therefore generalize  a Principal Component Analysis by extending this framework. This intermediary step was used, for example, by \cite{Bigot} to extend the definition of PCA on the space of probability measures on $\mathbb{R}$, and by \cite{Bigot_siam} to develop a suitable geometric PCA of a set of images.
A popular strategy to define moments in complex metric spaces is to use the variational interpretation of means (or barycenters), which leads to the introduction of Fr\'echet (or Karcher) means (see Section \ref{sec:frechet} for an accurate definition). This approach has been introduced in the seminal contribution \cite{frechet} that makes it possible to define $p$-means over any metric probability space.

 
The use  of Fr\'echet means has met with  great interest, especially in the field of bio-statistics and signal processing, although mathematical and statistical derivations around this notion constitute  a growing field of interest.
\begin{itemize}
\item[$\bullet$]
In  continuous domains, many authors have recently proposed limit theorems on the empirical  Fr\'echet sample mean (only a sample of size $n$ of the probability law is observed) towards its population counterpart. 
These works were mostly guided by applications to continuous manifolds that describe  shape spaces introduced in \cite{dryden}. For example, \cite{le} establishes the consistency of the population Fr\'echet mean and derives  applications in the Kendall space. The study of \cite{intrinsic} establishes the consistency of  Fr\'echet empirical means and derives its asymptotic distribution when $n \longrightarrow + \infty$, whereas the observations live in  more general Riemannian manifolds. Finer results can be obtained in some non-parametric restrictive situations  (see, \textit{e.g.}, \cite{bigot_gadat,bigot_gendre,bontemps_gadat} dedicated to the so-called \textit{shape invariant model}).
Many applications in  various domains involving signal processing can also be found: ECG curve analysis \cite{bigot2} and image analysis \cite{pennec,allasson}, to name a few.
\item[$\bullet$]
Recent works treat Fr\'echet means in a discrete setting, especially when dealing with 
 phylogenetic trees that have an important hierarchical structure property. In particular, \cite{le-owen} proposed a central limit theorem in this discrete case, whereas \cite{owen} used an idea of Sturm   for spaces with non-positive curvature to define an algorithm for the computation of the population Fr\'echet mean. Other works deal with the averaging of  discrete structure sequences such as diagrams using the Wasserstein metric (see, \textit{e.g.}, \cite{munch}) or graphs \cite{ginestet}.
\end{itemize}
Our work deviates from the above-mentioned point since we build an algorithm that recovers the Fr\'echet mean of a weighted graph while observing an infinite sequence of nodes of the graph, instead of finding the population Fr\'echet mean of a set of discrete structures, as proposed in \cite{le-owen,owen,munch,ginestet}.
 Our algorithm uses recent contributions on simulated annealing (\cite{miclo1,miclo2}). It relies on a continuous-time noisy simulated annealing Markov process on graphs, as well as a second process that accelerates 
 and homogenizes the updates of the noisy transitions in the simulated annealing procedure. 

The paper is organized as follows: Section \ref{sec:algo} highlights the different difficulties raised by the computation of the Fr\'echet mean of weighted graphs and
describes the algorithm we propose. In Section \ref{sec:quantum_diffusion}, theoretical backgrounds to understand the behavior of our method are developed on the basis of this algorithm, and Section \ref{sec:quantum_algo} states our main results.
Simulations and numerical insights are then given in Section \ref{sec:numeric}. The convergence of the algorithm is theoretically established in Section \ref{sec:section_proof},  whereas Section \ref{sec:functional} describes functional inequalities in quantum graphs that will be introduced below.

\section*{Acknowledgments.}
The authors gratefully acknowledge Laurent Miclo for his stimulating discussions and helpful comments throughout the development of this work, and Nathalie Villa-Vialaneix for her interest and advice concerning simulations. The authors are also indebted to \textit{zbMATH} for making their database available to produce numerical simulations. 

\section{Stochastic algorithm on quantum graphs}\label{sec:algo}

We consider $G=(V,E)$ a finite connected and undirected  graph  with no loop, where $V=\{1,\ldots,N\}$ refers to the $N$ vertices (also called nodes) of $G$, and $E$ the set of edges that connect some couples of vertices in $G$. 

\subsection{Undirected weighted graphs}

The structure of $G$ may be described by the adjacency matrix $W$ that gives a non-negative weight to each edge  $E$ (pair of connected vertices), so that $W=(w_{i,j})_{1 \leq i,j\leq N}$ while $w_{i,j}= + \infty$ if there is no direct link between node $i$ and node $j$. $W$ indicates the length of each direct link in $E$:  a small positive value of $w_{i,j}$ represents a  small length of the edge $\{i,j\}$. We assume that $G$ is undirected (so that the adjacency matrix $W$ is symmetric) and connected: for any couple of nodes in $V$, we can always find a path that connects these two nodes.
Finally, we assume that $G$ has no self loop. Hence, the matrix $W$ satisfies:
$$
\forall i\neq j \quad w_{i,j}=  w_{j,i} \qquad \text{and} \qquad \forall i \in V \quad w_{i,i}=0
$$
We define $d(x,y)$ as the geodesic distance between two points $(x,y) \in V^2$, which is the length of the shortest path between them. The length of this path is given by the addition of the length of traversed edges: 
$$
\forall (i,j) \in V^2 \qquad 
d(i,j) = \min_{i=i_1\rightarrow i_2 \rightarrow \ldots \rightarrow i_{k+1}=j} \sum_{\ell=1}^{k-1} w_{i_{\ell},i_{\ell+1}}.
$$

When the length of the edges is constant and equal to $1$, it simply corresponds to the number of traversed edges. Since the graph is connected and finite, we introduce the definition of the diameter of $G$:
$$
\mathcal{D}_G := \sup_{(x,y) \in V^2} d(x,y).
$$

To define a barycenter of a graph, it is necessary to introduce a discrete probability distribution $\nu$ over the set of vertices $V$. This probability distribution is used to measure the influence of each node on the graph.

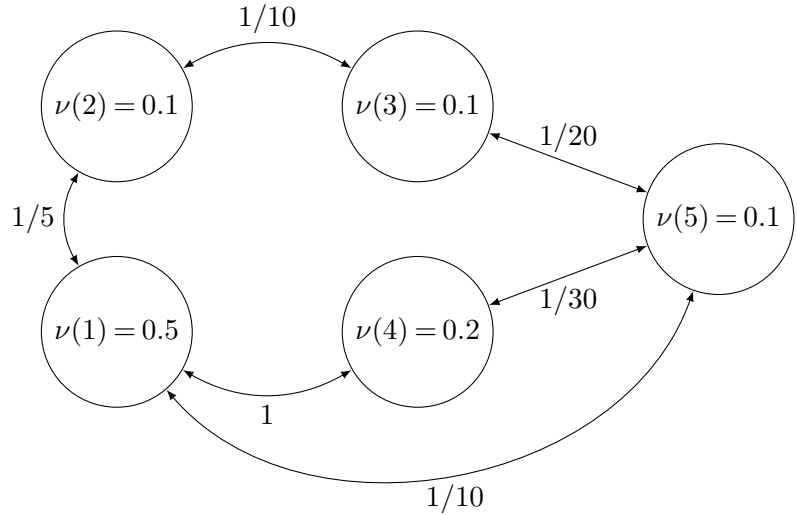
\begin{figure}
\begin{center}
\begin{tikzpicture}
    \node[shape=circle,draw=black] (A) at (0,0) { $\nu(1)=0.5$};
    \node[shape=circle,draw=black] (B) at (0,3) { $\nu(2)=0.1$};
    \node[shape=circle,draw=black] (C) at (4,3) {$\nu(3)=0.1$};
    \node[shape=circle,draw=black] (D) at (4,0) {$\nu(4)=0.2$};
    \node[shape=circle,draw=black] (E) at (8,1.5) {$\nu(5)=0.1$};

\draw[<->,>=latex] (A) edge[bend left] node[midway,left]{$1/5$} (B) ;
\draw[<->,>=latex] (B) edge[bend left] node[midway,above]{$1/10$} (C) ;
\draw[<->,>=latex] (A) edge[bend right] node[midway,below]{$1$} (D) ;
\draw[<->,>=latex] (A) edge[bend right=60] node[midway,below]{$1/10$} (E) ;
\draw[<->,>=latex] (E) edge[bend right=0] node[midway,above]{$1/20$} (C) ;
\draw[<->,>=latex] (E) edge[bend right=0] node[midway,below]{$1/30$} (D) ;

\end{tikzpicture}
\end{center}
\caption{\label{fig:ex_graphe}Example of a weighted graph between five authors, where  the probability mass  is $\{0.5,0.1,0.1,0.2,0.1\}$. In this graph, Author 1 shares five publications with Author 2 and ten publications with Author 5, and so on.}
\end{figure}

\begin{example}
Let us consider a simple scientometric example illustrated in Figure \ref{fig:ex_graphe} and consider a ``toy" co-authorship relation that could be obtained in a subgraph  of a collaboration network like \textit{zbMATH}\footnote{\url{https://zbmath.org/}}. If two authors $A$ and $B$ share $k_{A,B}$ joint papers, it is a reasonable choice to use a weight $w_{A,B} = \phi(k_{A,B})$, where $\phi$ is a convex function satisfying $\phi(0)=+ \infty,  \phi(1)=1$ and $\phi(+\infty)=0$. This means that no joint paper between $A$ and $B$ leads to the absence of a direct link between $A$ and $B$ on the graph. On the contrary, the more papers there are between $A$ and $B$, the closer $A$ and $B$ will be on the graph. Of course, this graph may be embedded in a probability space with the additional definition of a probability distribution over the authors that can be naturally proportional to the number of citations of each author. This is a generalization of the Erd\"os graph. Note that this type of example can  also be encountered when dealing with movies and actors, leading, for example, to the Bacon number and graph (see the website: \url{www.oracleofbacon.org/}).
\end{example}

\subsection{Fr\'echet mean of an undirected weighted graph}\label{sec:frechet}
Following the simple remark that the $p$-mean of any distribution $\nu$ of $\R^d$ is the point that minimizes
$$
x \longmapsto \mathbb{E}_{Z \sim \nu} [ |x-Z|^p],
$$
it is legitimate to be interested in the Fr\'echet mean of a graph $(G,\nu)$ where $\nu$ refers to the probability distribution over each node. 
The Fr\'echet mean is introduced in \cite{frechet} to generalize this variational approach to any metric space.
Although we have chosen to restrict our work to the case of $p=2$, which corresponds to the Fr\'echet mean definition, we believe that our work could be generalized to any value of $p \geq 1$.
 If $d$ denotes the geodesic distance w.r.t. $G$, we are interested in solving the following minimization problem:
\begin{equation}\label{eq:frechet_def}
M_{\nu} :=  \arg \min_{x \in E} U_{\nu}(x) \qquad \text{where} \qquad \um(x)=\frac{1}{2}\inte_{G} d^2(x,y)\nu(y).
\end{equation}
Hence,  the Fr\'echet mean of $(G,\nu)$, denoted $M_{\nu}$, is  the set of all possible minimizers of $\um$. Note that this set is not necessarily a singleton and this uniqueness property  generally requires some additional topological assumptions (see \cite{miclo1}, for example). At this point, we want to make three important remarks about the difficulty of this optimization problem:

\begin{itemize}
\item 
The problem of finding $M_{\nu}$ involves the minimization of $\um$, which is a non-convex function with the possibility of numerous local traps. To our knowledge, this problem cannot be efficiently solved using either a relaxed solution or using a greedy/dynamic programing algorithm (in the spirit of the Dijkstra method that makes it possible to compute geodesic paths \cite{dijkstra}).
\item  It is thus natural to think about the use of a global minimization procedure, and, in particular, the simulated annealing (S.A. for short) method. S.A. is a standard strategy to minimize a function over discrete spaces and its computational cost is generally high. It relies on an inhomogeneous Markov process that evolves on the graph with a transition kernel depending on the energy estimates $\um$. 
In other words, it requires the computation of $\um$  that depends on an integral w.r.t. $\nu$.
Since we plan to handle large graphs, this last dependency can be a very strong limitation.

\item In some cases, even the global knowledge of $\nu$ may not be realistic, and the importance of each node can only be revealed through i.i.d. sequential arrivals  of new observations in $E$ that are distributed according to $\nu$. This may be the case, for example, if we consider a probability distribution over $E$ that is cropped while gathering interactive forms on a website.
\end{itemize}

\subsection{Outline of Simulated Annealing (S.A.)}

The optimization with S.A. introduces a Markov random dynamical system that evolves either in continuous time (generally for continuous spaces) or in discrete time (for discrete spaces). When dealing with a discrete setting, S.A. is based on a  Markov kernel proposition $L(.,.): V \times V \longrightarrow [0,1]$ related to the $1$-neighborhoods of the Markov chain, as introduced in \cite{Hajeck}. It is based on an inhomogeneous Metropolis-Hastings scheme, which is recalled in Algorithm \ref{algo:MHSA}. We can derive asymptotic guarantees of the convergence towards a minimum of $U$ as soon as the cooling schedule is well chosen.

\SetInd{0.25em}{1em}
\begin{algorithm}[H]
 \KwData{Function $U$. Decreasing temperature sequence $(T_k)_{k \geq 0}$}

  \textbf{Initialization:} $X_0 \in V$\;
\For{$k = 0 \ldots + \infty$}{
\qquad Draw $x' \sim L(X_k,.)$ and compute $p_k = 1 \wedge \left\{ e^{T_k^{-1} [U(X_k)-U(x')]} \frac{L(x',X_k)}{L(X_k,x')}\right\}$

\qquad Update $X_{k+1}$ according to
$X_{k+1}=\begin{cases} x' \mbox{ with probability } p_k\\
X_k \mbox{ with probability } 1-p_k
\end{cases}$
 }
 
 \textbf{Output:} $\lim_{k\longrightarrow + \infty} X_k$
 \caption{M.-H. Simulated Annealing\label{algo:MHSA}}
\end{algorithm}

When dealing with a continuous setting,  S.A. uses a drifted diffusion with a vanishing variance $(\epsilon_t)_{t \geq 0}$ over $V$, or an increasing drift coefficient $(\beta_t)_{t \geq 0}$. We refer to \cite{stroock_sa,miclo_sa} for details and we  recall its Langevin formulation in Algorithm \ref{algo:LSA}.

\SetInd{0.25em}{1em}
\begin{algorithm}[H]
 \KwData{Function $U$. Increasing inverse temperature $(\beta_t)_{t \geq 0}$}

 \textbf{ Initialization:} $X_0 \in V$\;
 
 \qquad $
\forall t \geq 0 \qquad dX_t = - \beta_t \nabla U(X_t) dt + dB_t
$ 

\textbf{Output:} $\lim_{t\longrightarrow + \infty} X_t$
 \caption{Langevin Simulated Annealing\label{algo:LSA}}
\end{algorithm}

In both cases, we can see that S.A. with $U= \um $ given by \eqref{eq:frechet_def} involves the computation of the value of $U$ in Line 4 of Algorithm \ref{algo:MHSA}, or the computation of $\nabla U$ in Line 2 of Algorithm \ref{algo:LSA}. These two computations are problematic for our Fr\'echet mean problem: the integration over $\nu$ is intractable in the situation of large graphs and we are naturally driven to consider a noisy version of  S.A. 
A possible alternative method for this problem is to use a homogenization technique: replacing $\um$ in the definition of $p_k$ by $U_y(\cdot)=d^2(\cdot,y)$, \textit{where $y$ is a value from an i.i.d. sequence $(Y_n)_{n \geq 0}$ distributed according to $\nu$} .
Such methods have been developed in \cite{pflug} as a modification of Algorithm \ref{algo:MHSA} with  an additional Monte-Carlo step in Line 4, when $U_Y(x)$ follows a Gaussian distribution centered around the true value of $U(x) = \mathbb{E}_{Y \sim \nu} U_Y(x)$. This approach is still problematic in our case since the Gaussian assumption on the random variable $d^2(x,Y)$ is  unrealistic here. 
Another limitation of this MC step relies on the fact that it requires a batch average of several $(U(x,Y_j))_{1 \leq j \leq n_k}$ where $n_k$ is the number of observations involved at iteration $k$, although we also plan to develop an algorithm that may be adapted to on-line arrivals  of the observation $(Y_n)_{n \geq 0}$.
Lastly, it is important to observe that the non-linearity of the exponential prevents the use of only \textit{one} observation $Y_k$ in the acceptation/reject ratio involved in the S.A. since it does not lead to an unbiased evaluation of the true transition:
$$
\mathbb{E}_{Y \sim \nu} \left[ e^{T_k^{-1}[ U_{Y}(X_k) - U_Y(x')]} \right] \neq 
 e^{T_k^{-1} \mathbb{E}_{Y \sim \nu} [ U_{Y}(X_k) - U_Y(x')]}  \, .
$$

This difficulty does not arise in the homogenization of the simulated annealing algorithm in the continuous case since the exponential is replaced by a gradient, \textit{i.e.},  the process we use is a  Markov process of the form:
$$\de X_t=-\beta_t \nabla U_{Y_t}(X_t) dt + dB_t,$$ 
where $B_t$ is a ``Brownian motion" on the graph $G$,
$\beta_t^{-1}$ refers to the inverse of the temperature, and $Y_t$ is a continuous time  Markov process obtained from the sequence $(Y_n)_{n\in \mathbb{N}}$.
 This point motivates the introduction of the quantum graph induced by the initial graph. Of course, dealing with a continuous diffusion over a quantum graph $G$ deserves   special theoretical attention, which will be given in Section \ref{sec:quantum_diffusion}.

\subsection{Homogenized S.A. algorithm on a quantum graph}

We now present the proposed algorithm for estimating Fr\'echet means. To do so, we first introduce the quantum graph  $\Gamma_G$  derived from $G=(V,E)$ that corresponds to the set of points living inside the edges $e \in E$ of the initial graph. Once an orientation is arbitrarily fixed for each edge of $E$, the location of a point in $\Gamma_G$  depends on the choice of an edge $e \in E$ and on a coordinate $x_e \in [0,L_e]$ where $L_e$ is the length of  edge $e$ on the initial graph. The coordinate $0$ then refers to the initial point of $e$ and $L_e$ refers to the other extremity.

While considering the quantum graph $\Gamma_G$, it is still possible to define the geodesic distance between any point $x \in \Gamma_G$ and any node $y \in G$. In particular, when $x \in V$, we use the initial definition of the geodesic distance over the discrete graph, although when $x \in e \in E$ with a coordinate $x_e \in [0,L_e]$, the geodesic distance between $x$ and $y$ is:
$$
d(x,y) = \{ x_e + d(e(0),y)\} \wedge \{ L_e - x_e + d(e(L_e),y)\}  \, .
$$
This definition can be naturally generalized to any two points of $\Gamma_G$, enabling us to consider the metric space $(\Gamma_G,d)$.  

Consider a positive, continuous and increasing function  $t \longmapsto \alpha_t$ such that:

$$
\lim_{t \longmapsto + \infty} \alpha_t = + \infty \qquad \text{and} \qquad \forall t \geq 0 \quad  \beta_t =o(\alpha_t).
$$
We introduce $(N^{\alpha}_t)_{t \geq 0}$, an inhomogeneous Poisson process over $\mathbb{R}_+$  with intensity $\alpha$. It is standard to represent $N^{\alpha}$ through a  homogeneous Poisson process $H$ of intensity $1$ using the relationship:
$$ \forall t \geq 0 \qquad N_t^{\alpha}=H_{h(t)}, \mbox{ where } h(t)=\int_{0}^{t}\alpha_s\de s .$$

Using this accelerated process $N^\alpha$, our optimization algorithm over $\Gamma_G$ is based on the Markov process that solves the following stochastic differential equation over $\Gamma_G$:
\begin{equation}\label{eq:def_X}
\begin{cases} X_0 \in \Gamma_G\\
dX_t  = - \beta_t \nabla U_{Y_{N_t^{\alpha}}}(X_t) dt + dB_t.
\end{cases}
\end{equation}
Using the definition of $N_t$ and the basic properties of a Poisson Process, it can be observed that for all $\epsilon>0$ and all $t\ge 0$:
$$\mathbb{E}\left[\dfrac{N_{t+\epsilon}-N_t}{\epsilon}\right]=\frac{1}{\epsilon}\int_{t}^{t+\epsilon}\alpha_s\de s$$
Hence, $\alpha$ should be understood as the speed of new arrivals in the sequence $(Y_n)_{n\in\mathbb{N}}$.

Our definition \eqref{eq:def_X} is slightly inaccurate: for any $y \in V(G)$, the function $U_y$ is $\mathcal{C}^1$ except at a finite set of points of $\Gamma_G$:
\begin{itemize}
\item  This can be the case inside an edge $x \in e$ with $x_e \in (0,L_e)$ when at least two different geodesic paths from $x$ to $y$ start in opposite directions. Then, we define $\nabla U_y(x) = 0$. 
\item This can also be the case at a node $x \in V$ when several geodesic paths start from $x$ to $y$. In that case, we once again  arbitrarily impose a null value for the ``gradient" of $U_y$ at node $x$. This prior choice will not have any influence on the behavior of the algorithm.
\end{itemize}

We will first detail below the theoretical objects involved in \eqref{eq:def_X}, then we will describe an efficient discretization of \eqref{eq:def_X} that makes it possible to derive our practical optimization algorithm.

\SetInd{0.25em}{1em}
\begin{algorithm}[H]
 \KwData{Function $U$. Increasing inverse temperature $(\beta_t)_{t \geq 0}$. Intensity $(\alpha_t)_{t \geq 0}$}

  \textbf{Initialization:} Pick $X_0 \in \Gamma_G$. \;

$T_0=0$ \;
  
\For{$k = 0 \ldots + \infty$}{

\While{$N_t^{\alpha}=k$}{

\qquad $X_t$ evolves as a Brownian motion, relatively to the structure of $\Gamma_G$ initialized at $X_{T_k}^{-}$.
}

$T_{k+1} := \inf \{t : N_t^{\alpha} = k+1\}$\;

At time $t=T_{k+1}$, draw $Y_{k+1} = Y_{N_t^{\alpha}}$ according to $\nu$.

The process $X_t$ jumps from $X_t^-$ towards $Y_{k+1}$: 
\begin{equation}\label{eq:jumpX}X_t=X_t^{-} + \beta_t \alpha_t^{-1} \vec{X_tY_{N_t^{\alpha}}},\end{equation} where $\overrightarrow{X_tY_{N_t^{\alpha}}}$ represents the shortest (geodesic) path from $X_t$ to $Y_{N_t^{\alpha}}$ in $\Gamma_G$.
 }
 
 \textbf{Output:} $\lim_{t\longrightarrow + \infty} X_t$
 \caption{Homogenized Simulated Annealing over a quantum graph\label{algo:HLSA}}
\end{algorithm}

\begin{figure}
\begin{center}
\begin{tikzpicture}

    \node[shape=circle,draw=black] (A) at (0,0) {$V_1$};
    \node[shape=circle,draw=black] (B) at (0,3) {$V_2$};
    \node[shape=circle,draw=black] (C) at (4,3) {$V_3$};
    \node[shape=circle,draw=black] (D) at (4,0) {$V_4$};
    \node[shape=circle,draw=black] (E) at (8,1.5) {$V_5$};

\draw[-,>=latex] (A) edge  node[midway,left]{} (B) ;
\draw[-,>=latex] (B) edge node[midway,above]{} (C) ;
\draw[-,>=latex] (A) edge  node[midway,below]{} (D) ;
\draw[-,>=latex] (A) edge  node[midway,below]{} (E) ;
\draw[-,>=latex] (E) edge[bend right=0] node[midway,above]{} (C) ;
\draw[-,>=latex] (E) edge[bend right=0] node[midway,below]{} (D) ;
\draw[-,>=latex] (E) edge[bend right=0] node[midway,above]{} (B) ;

\draw (3.05,2.7) node[below]{\textcolor{red}{+}};
\draw (3,2.5) node[below]{\textcolor{red}{$X_{T_1}^-$}};

\draw (0,2.25) node[below]{\textcolor{red}{+}};
\draw (0,2) node[left]{\textcolor{red}{$X_{T_1}$}};

\draw (0,-0.5) node[left]{\textcolor{blue}{$Y_{T_1}$}};
\draw (1,4.7) node[left]{\textcolor{red}{\tiny Jump of size $\{ \beta_{_{T_1}} \alpha^{-1}_{_{T_1}}\} \times d(X_{T_1}^-,Y_{T_1})$}};

\draw (3.05,2.43)[draw=red] ..controls +(-4,3) and +(-4,3).. (0,2);

\draw (0,1.45) node[below]{\textcolor{red}{+}};
\draw (0,2.25) [draw=red]-- (0,1.45);
\draw (0,1) node[left]{\textcolor{red}{$X_{T_2}^{-}$}};
\draw (-0.2,1.5) node[left]{\textcolor{red}{\tiny Brownian displacement}};

\draw (8.3,1.5) node[right]{\textcolor{blue}{$Y_{T_2}$}};

\draw (3.05,0.85) node[below]{\textcolor{red}{+}};
\draw (3.5,0.85) node[left]{\textcolor{red}{$X_{T_2}$}};

\draw (0,1.15)[draw=red] ..controls +(0,0.5) and +(0,0.5).. (3.05,0.65);

\draw (5,1.5) node[left]{\textcolor{red}{\tiny Jump of size $\{ \beta_{_{T_2}} \alpha^{-1}_{_{T_2}}\} \times d(X_{T_2}^-,Y_{T_2})$}};
\end{tikzpicture}
\end{center}
\caption{\label{fig:ex_SA}Schematic evolution of the Homogenized S.A. described in Algorithm \ref{algo:HLSA} over the quantum graph. We  first observe a jump at time $T_1$ towards $Y_{T_1}=V_1$ and a Brownian motion on $\Gamma_G$ during $T_2-T_1$. A second node is then sampled according to $\nu$: here, $Y_{T_2}=V_5$ and a jump towards $Y_{T_2}$ occurs at time $T_2$.}
\end{figure}

Algorithm \ref{algo:HLSA} could be studied following the road map of \cite{miclo3}. Nevertheless, this  implies serious regularity difficulties on the densities and the  Markov semi-group involved. Hence, we have chosen to consider Algorithm \ref{algo:HLSA} as a natural Euler explicit discretization of our Markov evolution \eqref{eq:def_X}: for a large value of $k$, the average time needed to travel from $T_k$ to $T_{k+1}$ is approximately $\alpha_{T_k}^{-1} \longrightarrow 0$ as $k\longrightarrow +\infty$. On this short time interval, the drift term in \eqref{eq:def_X} is the gradient of the squared geodesic distance between $X_{T_k}$ and $Y_{N^{\alpha}_{T_k}}$, which is approximated by our vector $\overrightarrow{X_tY_{k}}$, multiplied by $\beta_{T_k}$, leading to \eqref{eq:jumpX}.
$X_t$ now  evolves as a Brownian motion over $\Gamma_G$ between two jump times and this evolution can be simulated with a Gaussian random variable using a (symmetric) random walk when the algorithm hits a node of $\Gamma_G$. Figure \ref{fig:ex_SA} proposes a schematic evolution of $(X_t)_{t \geq 0}$ over a simple graph $\Gamma_G$ with five nodes.
We will prove the following result.

\begin{theorem}\label{theo:conv}
A constant $c^\star(U_{\nu})$ exists such that if $\alpha_t = 1+t$ and $\beta_t = b \log(1+t)$ with $b>\{c^\star(U_{\nu})\}^{-1}$, then
$(X_t)_{t \geq 0}$ defined in \eqref{eq:def_X} converges a.s. to  $M_{\nu}$ defined in \eqref{eq:frechet_def}, 
when $t$ goes to infinity.
\end{theorem}

A more rigorous form of this result and more details about the process defined in \eqref{eq:def_X}  are presented in Section $\ref{sec:quantum_algo}$. The proof of Theorem \ref{theo:conv} is deferred to Section \ref{sec:section_proof}.

\section{Inhomogeneous Markov process over $\Gamma_G$ }
This section presents the theoretical background needed to define the Markov evolution \eqref{eq:def_X}.
\subsection{Diffusion processes on quantum graphs}\label{sec:quantum_diffusion}
We adopt here the convention introduced in \cite{freidlin_sheu} and fix for any edge $e\in E$ of length $L_e$ an orientation (and parametrization $s_e$). This means that $s_e(0)$ is one of the extremities of $e$ and $s_e(L_e)$ is the other one. By doing so, we have determined an orientation for $\Gamma_G$.  

\paragraph{Dynamical system inside one edge}

Following the parametrization of each edge, we can define the second order elliptic operator $\mathcal{L}_{e}$ as $\forall (x,y,t) \in e \times V(G) \times \mathbb{R}_+$:
\begin{equation}\label{eq:L_e}
  \mathcal{L}_{e} f (x,y,t) = - \beta_t \nabla_x U_{y}(x) + \frac{1}{2} \Delta_x f(x,y) + \alpha_t \int_{G} [f(x,y')-f(x,y)] d\nu(y'),
\end{equation}
which is associated with \eqref{eq:def_X} when $x \in e$: the $x$-component follows a standard diffusion drifted by $\nabla_x U_y(.)$ inside the edge $e$, although the $y$-component jumps over the nodes of the initial graph $G$ with a jump distribution $\nu$ and a rate $\alpha_t$. 

Since the drift term $\nabla U_y(.)$ is measurable w.r.t. the Lebesgue measure over $e$ and the second-order part of the operator is uniformly elliptic, $\mathcal{L}_e$ uniquely defines  (in the weak sense) a diffusion process up to the first time the process hits one of the extremities of $e$ (see, \textit{e.g.}, \cite{ikeda_watanabe}), which leads to a Feller Markov semi-group.

\paragraph{Dynamical system near one node}

We adopt the notation of \cite{FW} and write $e \sim v$ when a vertex $v \in  V$ is an extremity of an edge $e \in E$. 
For any function $f$ on $\Gamma_G$, at any point $v \in V$, we can define the directional derivative of $f$ with respect to an edge $e \sim v$ according to the parametrization of $e$. If $n_v$ denotes the number of edges $e$ such that $e \sim v$, we then obtain $n_v$ directional derivatives designated as $(d_ef(v))$:
$$
d_ef(v) = \displaystyle\begin{cases}
\displaystyle\lim_{h \longrightarrow 0^+} \frac{f(s_e(h))-f(s_e(0))}{h}  \qquad\qquad\,\,\,\,\, \text{if} \qquad s_e(0)=v\vspace{1em}\\
\displaystyle\lim_{h \longrightarrow 0^-} \frac{f(s_e(L_e+h))-f(s_e(L_e))}{h} \qquad \text{if} \qquad s_e(L_e)=v. \\
\end{cases}
$$

 It is shown in \cite{FW} that 
general dynamics over quantum graphs depend on a set of positive coefficients:
$$ \mathcal{A} := \left\{(a_v,(a_{e,v})_{e \sim v}) \in \mathbb{R}_+^{1+n_v}\,  \text{s.t.} \, a_v+\sum_{e \sim v} a_{e,v} >0 :  v \in V  \right\}.$$
There then exists a one-to-one correspondence between $\mathcal{A}$ and the set of all possible continuous Markov Feller processes on $\Gamma_G$.
More precisely, 
if the global generator $\mathcal{L}_t$  is defined as
$$
\forall (x,y) \in \Gamma_G \times V \qquad 
\mathcal{L}_tf(x,y) = \mathcal{L}_{e,t}(f)(x,y) \qquad \text{when} \qquad x \in e,
$$
while $f$ belongs to the domain:
\begin{equation}\label{eq:glue}
\mathcal{D}(\mathcal{L}) := \left\{ \forall y \in G \quad f(.,y) \in \mathcal{C}^{\infty}_b(\Gamma_G)  \, : 
\forall v \in V \qquad \sum_{e \sim v}  a_e d_ef(v,y) = 0\right\},
\end{equation}
then the martingale problem is well-posed (see \cite{FW,EK}).
In our setting, the \textit{gluing conditions} are defined through the following set  of coefficients $\mathcal{A}$:
$$
\forall v \in V \qquad 
a_v = 0  \qquad \text{and} \qquad \forall e \sim v  \qquad a_{e,v}  = \frac{1}{n_v}.
$$
The gluing conditions defined in \eqref{eq:glue} induce the following dynamics: when the $x$ component of $(X_t,Y_t)_{t \geq 0}$ hits an extremity $v \in V$ of an edge $e$, it is instantaneously reflected in one of the $n_v$ edges connected to $v$ (with a uniform probability distribution  over the connected edges) while spending no time on $v$. 
Using the uniform ellipticity of $\mathcal{L}_e$ and the measurability of the drift term, Theorem 2.1 of \cite{FW} can be adapted, providing the well-posedness of the Martingale problem associated with $(\mathcal{L},\mathcal{D}(\mathcal{L}))$, and the next preliminary result
can then be obtained.

\begin{theorem}
The operator $\mathcal{L}$ associated with the gluing conditions $\mathcal{A}$ generates a Feller Markov process on $\Gamma_G \times V \times \mathbb{R}_+$, with continuous sample paths on the $x$ component. This process is weakly unique and follows the S.D.E. \eqref{eq:def_X} on each $e \in E$.
\end{theorem}

For simplicity, $\Delta_x$ will refer to the Laplacian operator with respect to the $x$-coordinate on the quantum graph $\Gamma_G$, using our gluing conditions \eqref{eq:glue} and the formalism introduced in \cite{FW}.

\subsection{Convergence of the homogenized S.A. over $\Gamma_G$\label{sec:quantum_algo}}

 \medskip

As mentioned earlier, we  use a homogenization technique that involves an auxiliary sequence of random variables $(Y_n)_{n \geq 1}$, which are distributed according to $\nu$.  
More specifically, the stochastic process  $(X_t,Y_t)_{y \geq 0}$ described above is depicted by its inhomogeneous Markov generator, which can be split into three parts:
$$\mathcal{L}_tf(x,y)=\lu f(x,y) + \ld f(x,y).$$
In the equality above, $\lu$ is the part of the generator that acts on $Y_t$:
\begin{equation}\label{eq:deflu}\lu f(x,y)= \alpha_t\inte [f(x,y')-f(x,y)]\nu (\de y'),\end{equation}
describing the arrival of a new observation $Y \sim \nu$ with a rate $\alpha_t$ at time $t$.
Concerning   the action on the $x$ component, the generator is:

\begin{equation}\label{eq:defld}\ld f(x,y)=\dfrac{1}{2}\triangle_x f(x,y)-\beta_t<\triangledown_x U_y, \triangledown_x f(x,y)>.\end{equation}

Since the couple $(X_t,Y_t)_{t \geq 0}$ is Markov with a renewal of $Y$ with $\nu$, it can be immediately observed that the $y$ component is distributed at any time according to $\nu$. 
We introduce the notation $m_t$ to refer to the distribution of the couple $(X_t,Y_t)$ at time $t$, and we define $n_t$ as the marginal distribution of $X_t$.
In the following, we will also need to deal with the conditional distribution of $Y_t$ given the position $X_t$ in $\Gamma_G$. We will refer to this probability distribution as $m_t(y|x)$.
To sum up, we have:
\begin{equation}\label{eq:not_mesure}
\mathscr{L}(X_t,Y_t) = m_t \quad \text{with} \quad n_t(x) dx = \int_{V} m_t(x,y) \dy
\quad \text{and}\quad m_t(y\,\vert x) := \mathbb{P}[Y=y \, \vert X_t = x].
\end{equation}

A traditional method for establishing the convergence of S.A. towards the minimum of a function $U_{\nu}$ consists in studying the evolution of the law of $(X_t)_{t \geq 0}$ and, in particular, its close relationship with the Gibbs field $\mub$ with energy $U_{\nu}$ and inverse temperature $\beta_t$:
 \begin{equation}\label{eq:gibbs}\mub=\dfrac{1}{\z}\exp(-\beta_tU_{\nu}(x)),\end{equation}
 where $\z$ is the normalization factor, \textit{i.e.}, $\z = \inte e^{-\beta_t U_{\nu}(x)}\de x$.

Using a slowly decreasing temperature scheme $t \longmapsto \beta_t^{-1}$, it is expected that the process $(X_t)_{t \geq 0}$ evolves sufficiently fast (over the state space) in the ergodic sense so that its law remains close to $\mu_{\beta_t}$. 
In addition, the Laplace method on the sequence $(\mu_{\beta_t})_{t \geq 0}$ ensures that the measure is concentrated near the global minimum of $U_{\nu}$ (see, for example, the large deviation principle associated with $(\mu_{\beta})_{ \beta \to + \infty}$ in  \cite{FWFW}). 

Hence, a natural consequence of the convergence  ``$\mathscr{L}(X_t) \longrightarrow \mu_{\beta_t}$" and of the weak asymptotic concentration of $(\mu_{\beta_t})_{t \geq 0}$ around $M_{\nu}$ would be the almost sure convergence of the algorithm towards the Fr\'echet mean:

$$ \lim_{t \longrightarrow + \infty} P(X_t \in M_{\nu}) = 1.$$ We refer to \cite{stroock_sa} for further details. In particular, a strong requirement for this convergence can be considered through the relative entropy of $n_t$ (the law of $X_t$) with respect to $\mu_{\beta_t}$:
\begin{equation}\label{eq:entropie}
J_t := KL(n_t \vert \vert \mu_{\beta_t}) = \int_{\Gamma_G} \log \left[ \frac{n_t(x)}{\mu_{\beta_t}(x)}\right] d n_t(x)
\end{equation}
The function $\beta_t$ will be chosen as a $\mathcal{C}^1$ function of $\mathbb{R}_+$, and since $\mu_{\beta}$ is a strictly positive measure over $\Gamma_G$, it implies that $t \longmapsto \mu_{\beta_t}(x)^{-1}$ is $\mathcal{C}^1(\mathbb{R}_{+} \times \Gamma_G)$. Moreover, $(t,x) \longmapsto n_t(x)$ follows the backward Kolmogorov equation, which induces a  $\mathcal{C}^1(\mathbb{R}_+ \times \Gamma_G)$ function. Since the semi-group is uniformly elliptic on the $x$-component, we have: $$\forall t >0 \qquad  \forall x \in \Gamma_G \qquad  n_t(x) >0.$$
On the basis of  these arguments, we can deduce:

\begin{proposition} Assume that $t \longmapsto \beta_t$ is $\mathcal{C}^1$, then $(t,x) \longmapsto n_t(x)$ defines a positive $\mathcal{C}^1(\mathbb{R}_+ \times \Gamma_G)$ function and $t \longmapsto J_t$ is differentiable for any $t>0$.
\end{proposition}

\medskip
\noindent
 If we define:
  $$\alpha_t= \lambda (t+1) \mbox{ and } \beta_t =b \log(t+1),$$
with $b$  a constant smaller than $ c^\star(U_\nu)^{-1}$, where $c^\star(U_\nu)$ is the maximal depth of a well containing  a local but not global minimum of $\um$, defined in \eqref{eq:cstar}, then our main result can be stated as follows:
\begin{theorem}\label{theo:conv}
For any constant $\lambda>0$ such that $\alpha_t =\lambda (t+1)$ and $\beta_t = b \log(1+t)$ with $b > \{c^\star(U_\nu)^{-1}\}^{-1}$, then:
$$
\lim J_t = 0 \qquad \text{as} \qquad t \longrightarrow +\infty.$$

\end{theorem}
This ensures that the process $X_t$ will almost surely converge towards $M_{\nu}$ and, therefore, towards a global minimum of $\um$. 

\medskip
The idea of the proof is to obtain a differential inequality for $J_t$, which implies its convergence towards $0$.
 It is well known that the Gibbs measure $\mub$ is the unique invariant distribution of the stochastic process that evolves only on the $x$ component, whose Markov generator is given by:
\begin{equation}\label{eq:lcd_def}\lcd(f)(x) = \dfrac{1}{2} \Delta_x f -\beta_t \langle \nabla_x f, \nabla_x \um \rangle.\end{equation}
Therefore, a natural step of the proof will be to control the difference between $\ld$ and $\lcd$ and to use this difference to study the evolution of $m_t$ and $n_t$. 
It  can be seen that $\lcd$ may be written as an average action of the operator thanks to the linearity of the gradient operator:
$$
\lcd(f) = \frac{1}{2} \Delta_x f - \beta_t \mathbb{E}_{y \sim \nu} \langle \nabla_x f,\nabla_x U_y \rangle.
$$ 

When $X_t=x$, we know that $Y_t$ is distributed according to the distribution $m_t(y \vert x)$. Consequently, the average action of $\ld$ on the $x$ component is:
\begin{equation}\label{eq:ltd_def}
\ltd(f) = \frac{1}{2} \Delta_x f - \beta_t \mathbb{E}_{y \sim m_t(. \vert x)} \langle \nabla_x f,\nabla_x U_y \rangle =  \frac{1}{2} \Delta_x f - \beta_t \intev \langle \nabla_x f,\nabla_x U_y \rangle m_t(y \vert x) \dy,
\end{equation}
whose expression may be close to that of $\lcd$ if $m_t(.\vert x)$ is close to $\nu$.

Thus, another important step is to choose appropriate values for $\alpha_t$ and $\beta_t$, \textit{i.e.}, to find the balance between the increasing intensity of the Poisson process and the decreasing temperature schedule, in order to quantify the distance between $\nu$ and  $m_t(.\vert x)$.
The main core of the proof brings together these two aspects and is detailed in Section \ref{sec:section_proof}.

\medskip

Another important step will be the use of functional inequalities (Poincar\'e and log-Sobolev inequalities) over $\Gamma_G$ for the measure $\mu_{\beta}$ when $\beta=0$ and $\beta \longrightarrow +\infty$. The proof of these technical results are given in Section $\ref{sec:sobolev_section}$.

\begin{corollary}
Assume that $\beta_t = b \log(t+1)$ with $b < c^\star(U_\nu)^{-1}$ and $\alpha_t = \lambda (t+1)^\gamma$ with $\gamma \geq 1$, then for any neighborhood $\mathcal{N}$ of $M_{\nu}$:
$$
\lim_{t \longrightarrow + \infty} \mathbb{P}[X_t \in \mathcal{N}] = 1.
$$
\end{corollary}
  
\underline{\textit{Proof:}}
The argument follows  from Theorem \ref{theo:conv}. Consider any neighborhood $\mathcal{N}$ of $M_{\nu}$. The continuity of $U_{\nu}$ shows that:
$$
\exists \delta >0 \quad \mathcal{N}^c \subset \{x \in \Gamma_G \,: U(x) > \min U + \delta \}.
$$
Hence, 
\begin{eqnarray*}
\mathbb{P}[X_t \in \mathcal{N}^c] &\leq& \mathbb{P}[U(X_t) > \min U +\delta]\\
 & = & \integ \mathbf{1}_{U(x) > \min U + \delta} \de n_t(x)\\
 & = & \integ \mathbf{1}_{U(x) > \min U + \delta} \de \mub(x) + \integ \mathbf{1}_{U(x) > \min U + \delta} [n_t(x) -  \mub(x)] \dx\\
 & \leq & \mub \left\{ U > \min U +\delta \right\} + 2 d_{TV}(n_t,\mub)\\
 & \leq& \mub \left\{ U > \min U +\delta \right\} + \sqrt{2 J_t},
\end{eqnarray*}
where we used the variational formulation of the total variation distance and the Csisz\'ar-Kullback inequality. As soon as $\lim_{t \longrightarrow +\infty} J_t = 0$, we can conclude the proof observing that $\mub \left\{ U > \min U +\delta \right\} \longrightarrow 0$ as $\beta_t \longrightarrow + \infty$.
\hfill $\square$

It can actually be proven (not shown here) that $lim_{t \longrightarrow + \infty} \mathbb{P}[X_t \in \mathcal{N}] = 1$ if and only if the constant $b$ is chosen to be lower than $c^\star(U_\nu)^{-1}$. This means that this algorithm does not allow faster cooling schedules than the classical S.A. algorithm. Nevertheless, this is a positive result since this homogenized S.A. can be numerically computed quickly and easily on large graphs. Finally, we should  consider this result to be theoretical. However, in practice, the simulation of this homogenized S.A. is performed during a finite horizon time and efficient implementations certainly deserve a specific theoretical study following the works of \cite{Catoni} and \cite{Trouvé}.

\section{Numerical simulation results \label{sec:numeric}}

We now present different practical aspects of our strategy. After giving numerical details about how to make it usable on large graphs, plus different insights into its calibration, we will present results obtained on social network subgraphs of \textit{Facebook}\footnote{Subgraphs obtained on the Stanford Large Network Dataset Collection: \url{https://snap.stanford.edu/data/}} and a citation subgraph of \textit{zbMATH}.


\subsection{Algorithmic complexity}\label{ssec:AlgoComp}
It is widely recognized that the algorithmic complexity of numerical strategies dealing with graphs can be an issue. 
The number of vertices and edges of real-life graphs can  indeed be quite large. For instance, the \textit{zbMATH} subgraph of Section~\ref{sec:zbMATH} has 13000 nodes and approximately 48000 edges. In this context, it is worth justifying that the motion of $(X_t)_{t \geq 0}$ on the graph $\Gamma_G$ across the iterations of Algorithm~\ref{algo:HLSA} is reasonably demanding in terms of computational resources. We recall that this motion is driven by a Brownian motion (lines 4 to 6 of Algorithm~\ref{algo:HLSA}) and an attraction towards the  vertex $Y_{N_t^{\alpha}}$ (line 9 of Algorithm~\ref{algo:HLSA}) at random times $(T_k)_{k \geq 0}$. The number of vertices and undirected edges with non-null weights in $\Gamma_G$ is also $N$ and $|E|$, respectively. Note finally that $N-1 \leq |E|  \leq N (N-1)/2$ if $\Gamma_G$ has a unique connected component.

\medskip

\paragraph{Neighborhood structure}

A first computational issue that can arise when moving $X_t$ is to find all possible neighbors of a specific vertex $v$. This is indeed performed every time $X_t$ moves from one edge to another and directly depends on how $\Gamma_G$ is encoded in the memory.
Encoding  $\Gamma_G$ in a list of edges with non-null weights is common practice. In that case, the computer checks the vertex pairs linked by all edges to find those containing $v$, so the algorithmic cost is $2 |E|$. A second classic strategy is to encode the graph in a connectivity matrix. In this case, the computer has to go through all the $N$ indexes of the columns representing $v$ to find the non-null weights. We instead sparsely encode the graph in a list of lists: the main list has a size $N$ and each of its elements lists the neighbors of a specific vertex $v$. 
If the graph only contains edges with strictly positive weights, this strategy has a computational cost $N$, and in all of the other cases, the cost is $<N$. On average, the computational cost is $2 |E|/N$  so that this strategy is particularly advantageous for sparse graphs, where $|E|<<N (N-1)/2$,  which are common for the targeted applications. For instance,  $|E|=4*10^3$ and $N(N-1)/2=1.24*10^5$ on the smallest \textit{Facebook} subgraph of Section~\ref{sec:facebook}, and $|E|=4.8*10^4$ and $N(N-1)/2=8.44*10^7$ on the \textit{zbMATH} subgraph of Section~\ref{sec:zbMATH}. 

\medskip

\paragraph{Geodesic paths}

Another potential issue with our strategy is that it seeks at least one optimal path between the vertices  $X_t$ and $Y_{N_t^{\alpha}}$ at each jump time (line 9 of  Algorithm~\ref{algo:HLSA}). This may be efficiently done using a fast marching propagation algorithm (see, \textit{e.g.}, \cite{dijkstra} for details), where:
\begin{itemize}
\item[(i)]
the distance to $X_t$ is iteratively propagated on the whole graph $\Gamma_G$ until no more optimal distance to reach $X_t$ is updated, and 
\item[(ii)] considering the shortest path between $X_t$ and $Y_{N_t^{\alpha}}$.
\end{itemize}
In this case, step (i) is particularly time-consuming and \textit{cannot be reasonably performed at each step of Algorithm \ref{algo:HLSA}}. Fortunately, the algorithmic cost to compute the distance between all pairs of vertices is equivalent to that of computing step (i).
We then compute these distances once and for all at the beginning of the computations and store the result in a $N \times N$ matrix in the RAM of the computer \footnote{In our experiments, we used the \textit{all-shortest-paths} function of the Python library NetworkX.}.
We can therefore very quickly use these results at each iteration of the algorithm and, in particular, deduce the useful part of the geodesic paths involved to travel from $X_t$ to $Y_{N_t^{\alpha}}$.
The only limitation of this  strategy is that the distance matrix can be memory-consuming. 
It can therefore be used on small to large graphs but not on huge graphs  (typically when $N>10^5$) on current desktops. An extension to deal with this scalability issue is a current subject of research, and we will therefore not describe applications of Algorithm \ref{algo:HLSA} on huge graphs in this work.

\medskip

\subsection{Parameter tuning}\label{ssec:paramtunig}

Several parameters influence the behavior of the simulated process $X_t$. Some of them are directly introduced in the theoretical construction of the algorithm, \textit{i.e.}, the intensity of the Poisson process or the temperature schedule. Other ones come from the practical implementation of the algorithm, \textit{i.e.}, the maximal time up to which we generate $X_t$.
The theoretical result given in Theorem \ref{theo:conv}  gives an upper bound for the probability of $X_t$ to be not too distant from the set of global minima. 
This bound depends on  $\beta_t$ and $\alpha_t$ as well as on different characteristics of the graph such as its diameter and number of nodes. 
We then propose  an empirical strategy for parameter tuning, which we will use later in Section \ref{ssec:results}. 

For a graph $G$ with $N$ nodes and a diameter  $\mathcal{D}_G$, we define:
 $$T^*_{\max}=100+0.1 N \qquad \text{and} \qquad   \beta_t^*=2\log(t+1) / \mathcal{D}_G.$$
We then choose the intensity of the Poisson process so that it will generate a reasonably large amount of $Y_n$ at the end of the algorithm, depending on the discrete probability distribution $\nu$. More specifically, let $S^*$  be the average  number of jumping times between $T^*_{max}-1$ and $T^*_{max}$. We then set $S^*=1000$, which can be obtained by choosing:
$$\alpha_t^*=\lambda(2t+2) \qquad \text{with} \qquad \lambda=\frac{S^*}{2T^*_{\max}+1}.$$ 
\medskip

\subsection{Results}\label{ssec:results}


This section presents results obtained on graphs of different sizes and using different parameters. We used the empirical method of Section \ref{ssec:paramtunig} to define default parameters and altered them to quantify the sensitivity of our algorithm to parameter variations. 

Since the process $X_t$ lives on a quantum graph, its location at time $t$ is between two vertices, on an edge of $\Gamma_G$. 
However, our primary interest is to study the properties of the initial discrete graph and therefore, the output of the algorithm will be the vertex considered as the graph barycenter. To achieve this, we associate a frequency to each vertex. For a vertex $v$, this frequency is the portion of time during which $v$ was the closest vertex to the simulated process $X_t$. In our results, we consider frequencies computed over the last 10$\%$ of the iterations of the algorithm.

\subsubsection{Facebook subgraphs}\label{sec:facebook}

\paragraph{Experimental protocol}

We first tested our algorithm on three subgraphs of Facebook from the Stanford Large Network Dataset Collection: 
\textbf{(FB500)} has 500 nodes and 4337 edges and contains two obvious clusters; \textbf{(FB2000)} has 2000 nodes and 37645 edges and fully contains (FB500); and \textbf{(FB4000)} has 4039 nodes and 88234 edges and fully contains (FB2000).
%
For each of these subgraphs, we considered the probability measure $\nu$ as the uniform distribution over the graph's vertices and used a length of 1 for all edges. We also  explicitly computed the barycenter of these graphs using an exhaustive search procedure. We are able to do that because we considered a simplified case where the distribution over the nodes is uniform. For example, this exhaustive search procedure required approximately 6 hours for the \textbf{(FB4000)} subgraph. Nevertheless, it allowed us to compare our strategy with ground-truth results.

We used the strategy of Section \ref{ssec:paramtunig} to define default parameters adapted to each subgraph. We also tested different values for parameters $\beta$, $S$ and $T_{\max}$ in order to quantify their influence. We  repeated our algorithm 100 times for each parameter set to evaluate the algorithm stability.

In the tables representing quantitative results, \textit{Error} represents the number of times, out of 100, that the algorithm converged to a node different from the ground-truth barycenter. It is a rough indicator of the ability of the algorithm to locate the barycenter of the graph, that could be replaced by a measure of the average distance between the last iterations of the algorithm and the ground-truth barycenter (not shown in this work).
For each subgraph and parameter set,  column \textit{Av. time} contains the average times in seconds for the barycenter estimation, keeping in mind that the Dijkstra algorithm was performed once for all before the 100 estimations. This preliminary computation requires approximately 1, 30 and 80 seconds on a standard laptop.



\paragraph{Effect of $\beta$}

From a theoretical point of view, $(\beta_t)_{t \geq 0}$ should be chosen in relation to the constant $c^\star(U_{\nu})$, which is unknown in practice. Therefore, the practical choice of $(\beta_t)_{t \geq 0}$ is a real issue to obtain a good behavior of the algorithm.
Table \ref{tab:FB_beta} gives the results obtained on our algorithm  with different values of $\beta_t$. 

\begin{table}[htb]
 \small

\begin{center}
\begin{tabular}{|l||c|c|c||c|c|c||c|c|c|}
   \hline
                 & \multicolumn{3}{|c||}{FB500}       & \multicolumn{3}{|c||}{FB2000}      &\multicolumn{3}{c|}{FB4000} \tabularnewline  
    \hline
    $\beta$ & Error & Med. Freq. &  Av. time & Error & Med. Freq. &  Av. time &Error & Med. Freq. &  Av. time \tabularnewline
    \hline
    $\frac{1}{4}\beta^* $ & 15 $\%$ & 0.6042 & 0.95 s & 28 $\%$  & 0.3519 & 10.12 s & 27 $\%$ &0.3314 & 12.41 s
    \tabularnewline
    \hline
    $\frac{1}{2}\beta^* $ & 2 $\%$ &  0.4184 &  1.38 s & 1 $\%$ &    0.7268&  4.25 s & 5 $\%$ &   0.6534 &  14.11 s  \tabularnewline
    \hline
    $\beta^* $ & 0 $\%$ & 0.8008 &  1.48 s &0 $\%$ & 0.9418 & 12.36 s  & 1 $\%$ &0.8913 & 13.55 s \tabularnewline
    \hline
    $2\beta^* $ & 0 $\%$ & 0.8321 & 11.31 s & 0 $\%$ & 0.9892 & 8.52 s &0 $\%$ & 0.9647 & 16.79 s  \tabularnewline
    \hline
    $4\beta^* $ & 0 $\%$ & 0.8233 & 13.19 s &0 $\%$ & 0.9930 & 15.58 s& 0 $ \%$  & 0.9824 & 31.43 s \tabularnewline
    \hline
    $8\beta^* $ & 0 $ \%$ & 0.7717 & 2.36 s &0 $\%$ & 0.9750 & 23.25 s& 0 $ \%$  &0.9445 &44.46 s \tabularnewline
    \hline
    \end{tabular} 
\end{center}
\caption{
Experiments on the three Facebook subgraphs, while varying the values of $(\beta_t)_{t \geq 0}$.
\label{tab:FB_beta} }
\end{table}

\normalsize

We can observe that when $(\beta_t)_{t \geq 0}$ is too small, then the behavior of the algorithm  is deteriorated, revealing the tendancy of the process $(X_t)_{0 \leq t \leq T^*_{\max}}$ to have an excessively slow convergence rate towards its local attractor in the graph. Roughly speaking, in such a situation,  the process does not learn fast enough. 
When the value of $\beta$ is chosen in the range $[\beta^*;2\beta^*]$, we can observe a really good behavior of the algorithm: it almost always locates the good barycenter in a quite reasonable time of computation (less then 20 seconds for the largest graph). 
Finally, we can observe in the column, Med. Freq., that in most of the last iterations of the algorithm (more than 80$\%$), the process evolves around its estimated barycenter, so that the decision to produce an estimator is quite easy when looking at an execution of the algorithm.

\paragraph{Effect of $S$ and $T_{\max}$}

Table \ref{tab:FB_ST} gives the results obtained with our algorithm while using different values of $S$ and $T_{\max}$. 
As expected, we observe that increasing the ending time of simulation always improves the convergence rate (column Error in Table \ref{tab:FB_ST}) of the algorithm towards the right node. The behavior of the algorithm is also improved by increasing the value of $S$, which quantifies the number of arrivals of nodes observed along the averaging procedure. Of course, the counterpart of increasing both $S$ and $T_{max}$ is  an increasing cost of simulation (see column Av. time).

\begin{table}[htb]
  \small
\begin{center}
\begin{tabular}{|l|l||c|c|c||c|c|c||c|c|c|}
   \hline
              &   & \multicolumn{3}{c||}{FB500}       & \multicolumn{3}{|c||}{FB2000}      &\multicolumn{3}{c|}{FB4000} \tabularnewline  
    \hline
    $S$& $T_{\max}$ & Error   & Med. Freq. &  Av. time &Error   & Med. Freq. &  Av. time &Error   & Med. Freq. &  Av. time \tabularnewline
    \hline
$\frac{1}{2}S^*$ &$T_{\max}^*$ & 2 $\%$  & 0.7667 & 0.59 s   & 0 $\%$ &  0.9344 & 3.22 s & 0 $\%$ & $0.8614$ & 6.70 s\tabularnewline
    \hline
    $S^*$ &$2T_{\max}^*$ & 0 $\%$ &  0.8049 & 2.79 s & 0 $\%$ & 0.9610 & 9.30 s & 0 $\%$ & 0.8970 & 24.42 s\tabularnewline
    \hline
     $S^*$ &$4T_{\max}^*$ & 0 $\%$ & 0.8101 & 5.60 s &0 $\%$ &0.9677 & 22.21 s  & 0 $\%$ &0.9222 & 54.73 s \tabularnewline
    \hline
    $2S^*$ & $T_{\max}^*$ & 1 $\%$ & 0.8345 & 2.26 s  & 0 $\%$ & 0.9512 & 8.39 s & 0 $\%$ & 0.9062 & 20.95 s \tabularnewline
    \hline
    $2S^*$& $2T_{\max}^*$ & 0 $\%$ & 0.8361  & 4.57  s & 0 $\%$ & 0.9586  & 23.30 s & 0 $\%$ & 0.9121 & 23.30 s\tabularnewline
    \hline
     $2S^*$& $4T_{\max}^*$& 0 $\%$  & 0.8423 & 11.16 s & 0 $\%$  & 0.9735 & 11.16 s & 0 $\%$ & 0.9366 & 96.84 s\tabularnewline  
    \hline
    \end{tabular} 
\end{center}
\caption{
Experiments on the three Facebook subgraphs, while varying the values of $S$ and $T_{\max}$.
\label{tab:FB_ST}
 }
\end{table}   

\begin{figure}[htb!]
\begin{center}
\includegraphics[width=0.70\linewidth]{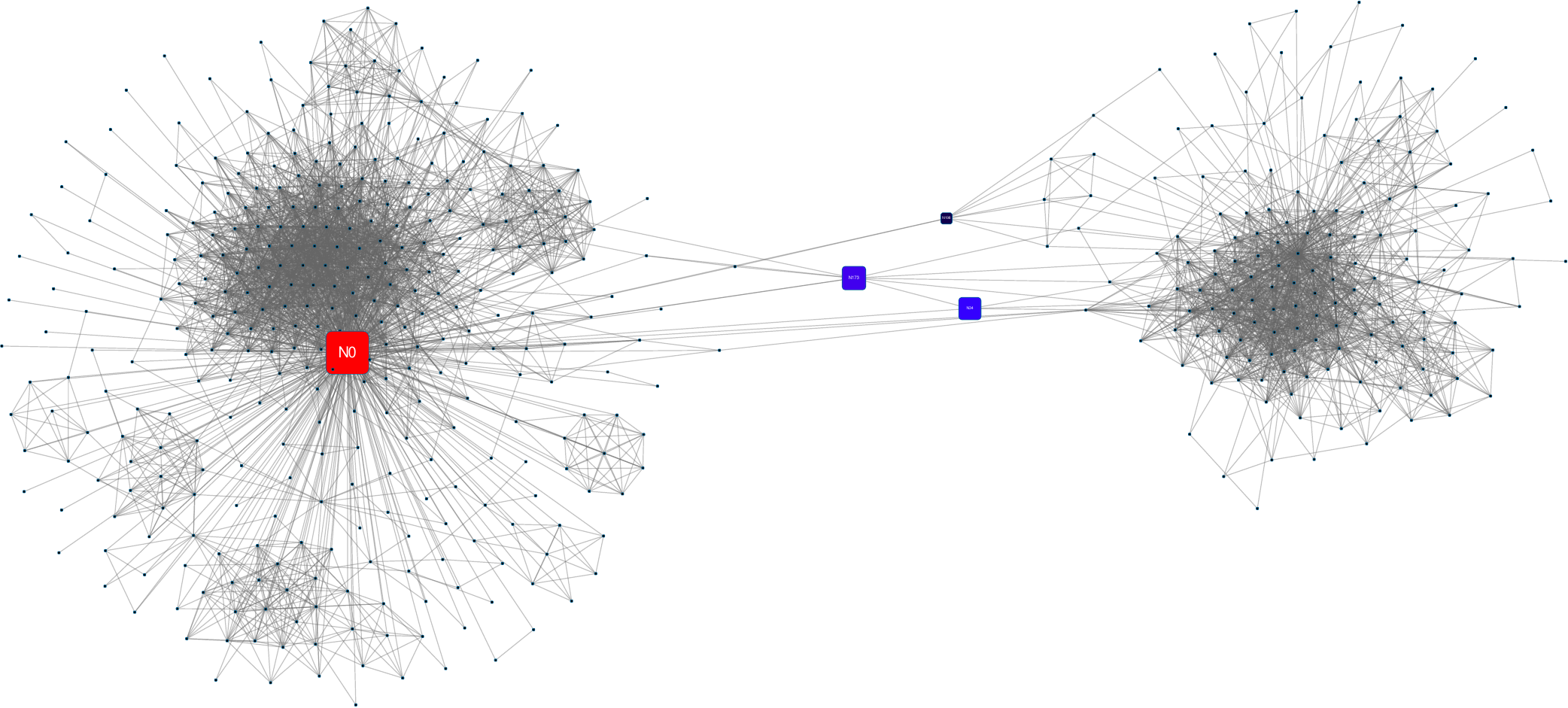}
\caption{
Region of interest presenting the results obtained on the FB500 graph.
}
\label{fig:pic_fb500}
\end{center}
\end{figure}

\normalsize
As an illustration of the (small) complexity of the \textit{Facebook} sub-graphs used for benchmarking our algorithm,
we provide a representation of the FB500 graph in Figure \ref{fig:pic_fb500}.
 This representation has been obtained with the help of \textit{Cytoscape} software and is not a result of our own algorithm. In Figure \ref{fig:pic_fb500}, the red node is the estimated barycenter, which is also the ground-truth barycenter located by a direct exhaustive computation.  The blue nodes are the ``second rank" nodes visited by our method.

\subsubsection{zbMATH subgraph}\label{sec:zbMATH}

The zbMATH subgraph built from zbMATH\footnote{\url{https://zbmath.org/authors/}} has been obtained by an iterative exploration of the co-authorship relationship in an alphabetical order. We thus naturally obtained a connected graph. This graph has been weighted by the complete number of citations obtained by each author, leading to the distribution probability $\nu\propto \sharp (\text{number of citations})$. Finally, all the edges in the graph are fixed to have a length of $1$.
This exploration was initialized on the entry of the first author's name (\textit{i.e.}, S. Gadat) and we stopped the process when we obtained 13000  nodes (authors) on the graph. This stopping criterion in the exploration of the zbMATH database corresponds to a technical limitation of $40$ GB memory required by the distance matrix obtained with the Dijsktra algorithm. In particular, this limitation and the starting point of the exploration induce an important bias in the community of authors used to build the subgraph from the zbMATH dataset: the researchers obtained in the subgraph are generally French and applied mathematicians. The graph is more or less focused on the following research themes: probability, statistics and P.D.E.
Consequently, the results provided below should be understood as an illustration of our algorithm and not as a bibliometric study!

\begin{figure}[htb!]
\begin{center}
\includegraphics[width=0.2\linewidth]{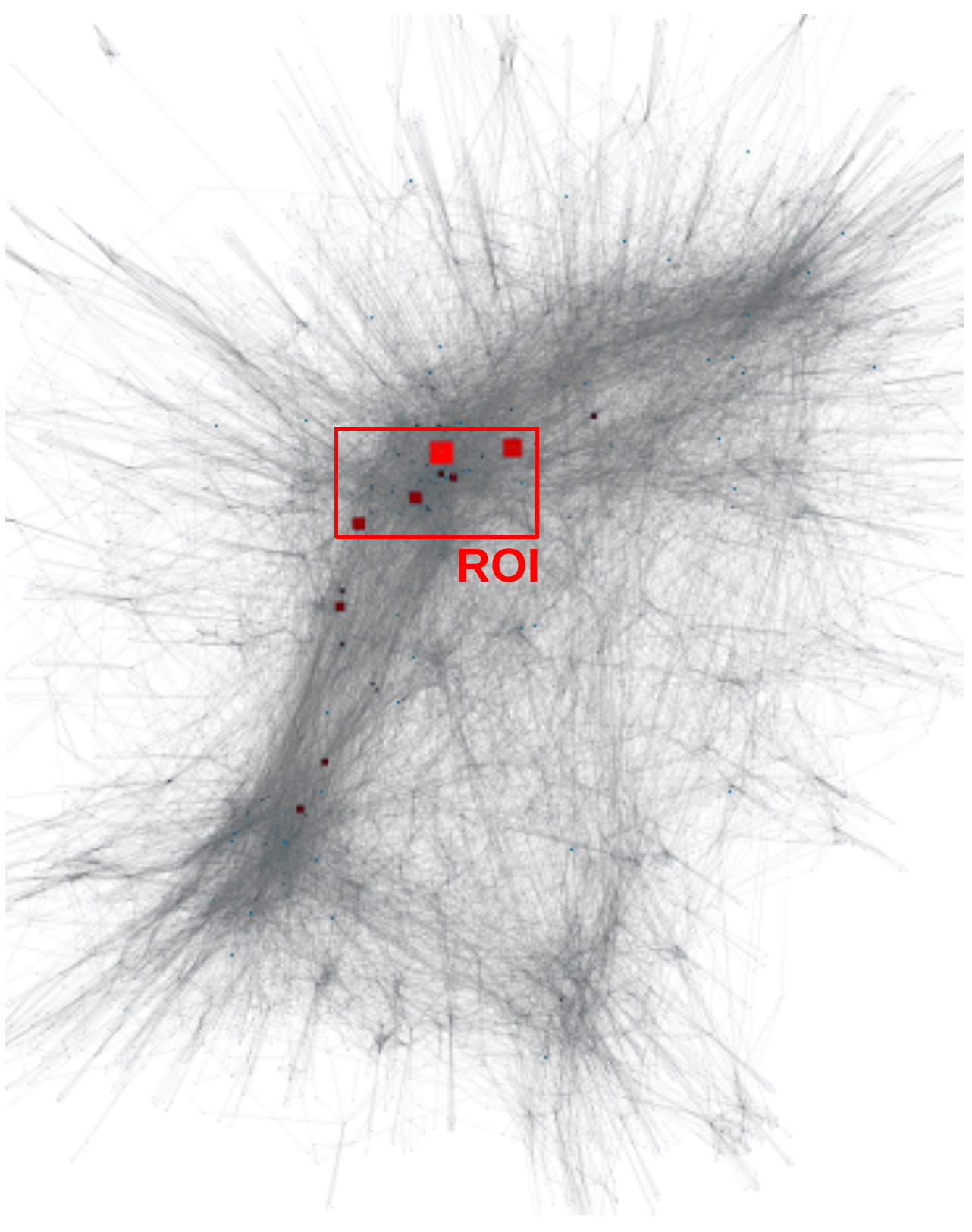} 
\includegraphics[width=0.75\linewidth]{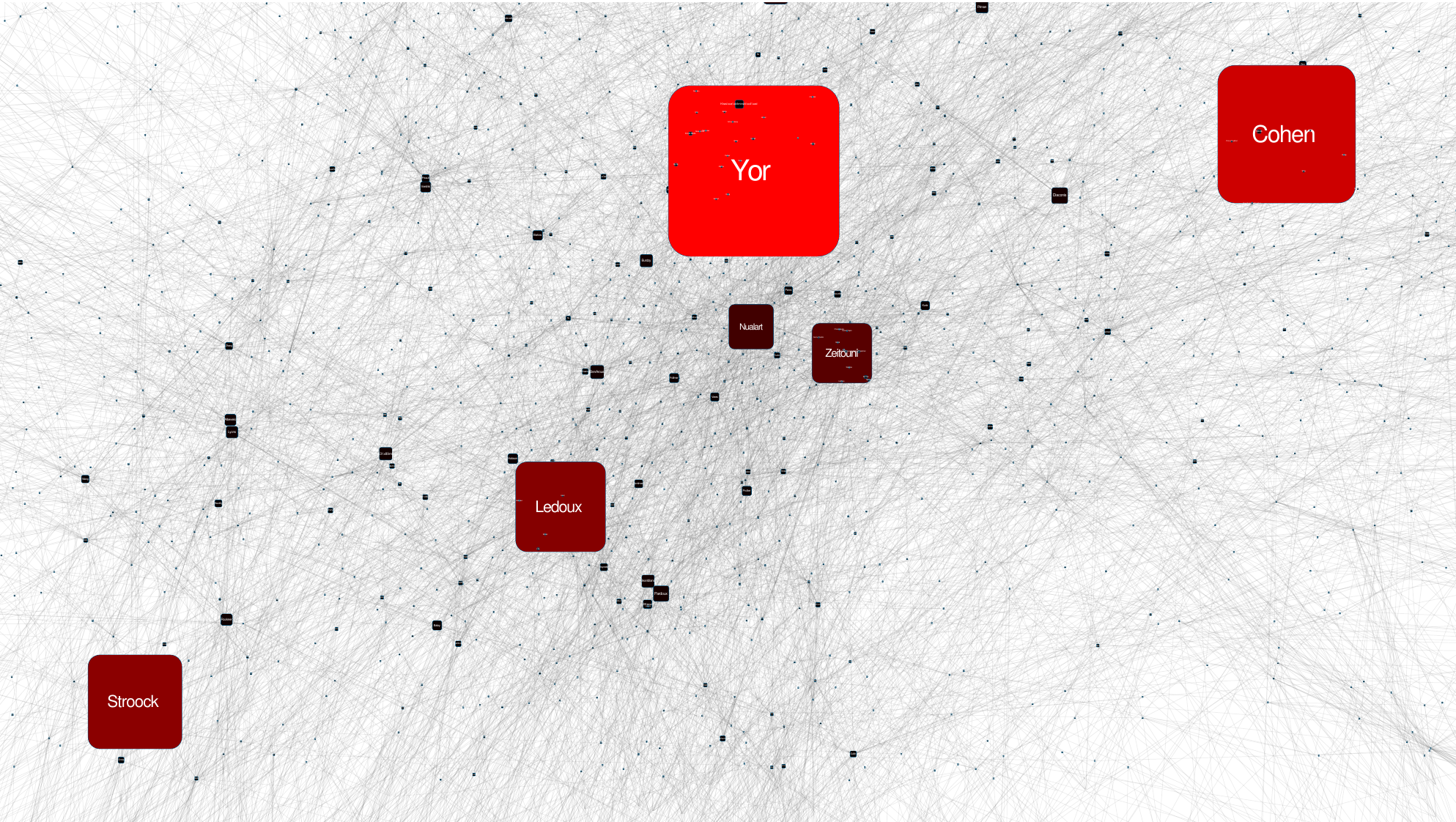} 
\caption{Left: General overview of the zbMath subgraph extracted for our experiments, 
containing approximately 13000 nodes. Right:
Region of interest presenting the main results obtained on the zbMATH subgraph.
}
\label{fig:zbmath}
\end{center}
\end{figure}

Our experiments rely on the same choice of parameters $(\beta_t^*)_{t \geq 0}, S^*$, $T^*_{max}$ and $\alpha_t^*$ as above, indicated at the beginning of Section \ref{ssec:paramtunig}. Again, our algorithm produces the same outputs but of course, in this situation, we do not know the ground-truth barycenter of the zbMath subgraph. 
In Figure \ref{fig:zbmath}, we present a representation of the subgraph obtained with \textit{Cytoscape} software, and a zoom on a region of interest (ROI for short) in Figure \ref{fig:zbmath}.

Again, the main nodes visited by our algorithm are represented with a  red square and the size of the used square is larger when the node is frequently visited. According to the results obtained on the \textit{Facebook} subgraphs, we then assume that the larger red square is the barycenter of the zbMath subgraph.

\begin{figure}[htb!]
\begin{center}

\includegraphics[width=0.75\linewidth]{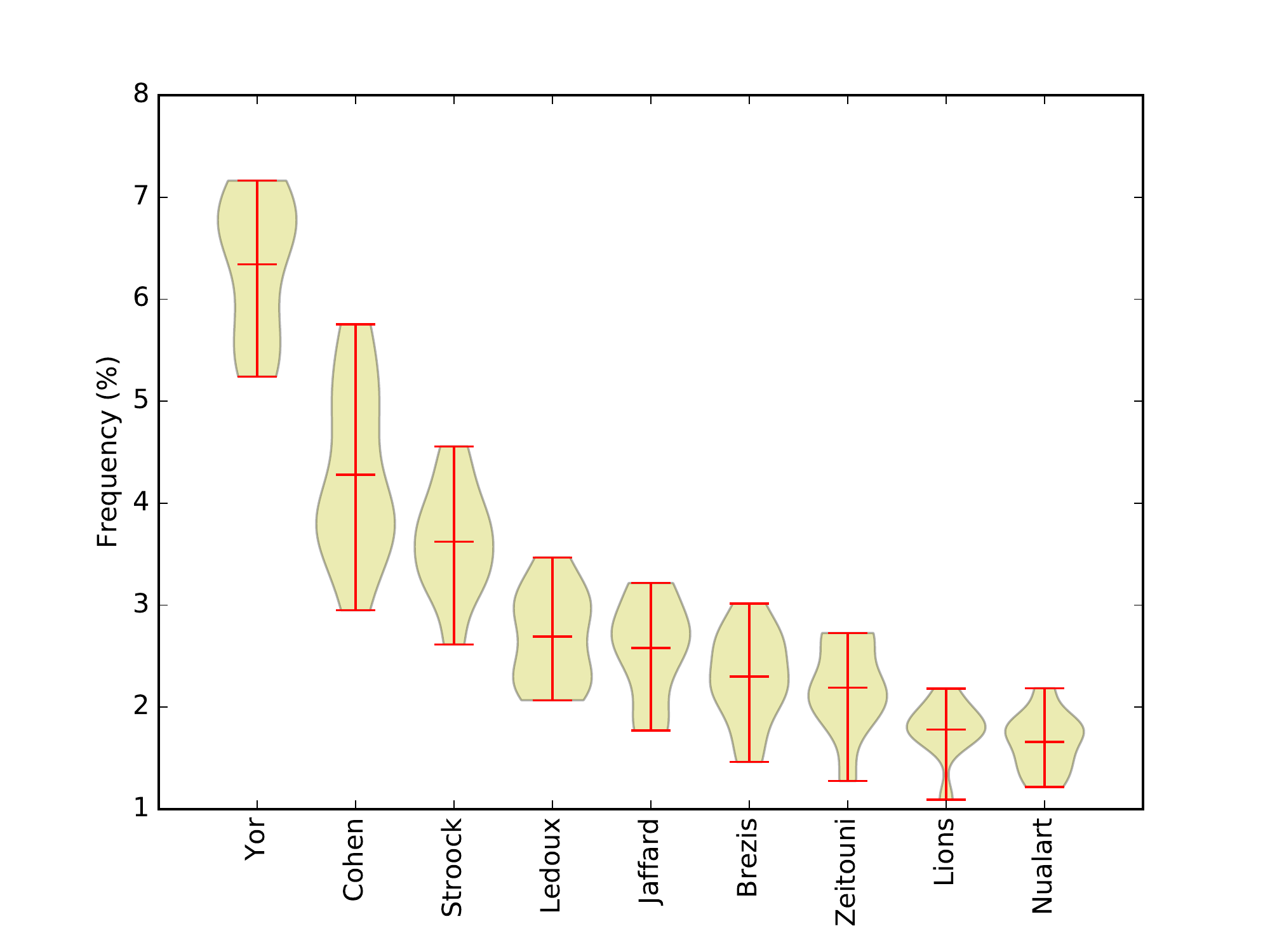} 
\caption{``Violin" plot of the occupation measure of the algorithm on the zbMATH subgraph for $10$ MC replications. The average frequency is located in the middle of the ``violin" plot, although the minimal and maximal values are shown in the extremity of the representation.
}
\label{fig:zbmath_violin}
\end{center}
\end{figure}

It is also possible to assert the robustness of our method with respect to several Monte Carlo runs of our algorithm. We have produced some boxplots for each of the main authors located in the subgraph, according to the occupation measure of the process over the last $10\%$ of the iterations with $10$ Monte-Carlo replications. These ``violin" plots are represented in Figure \ref{fig:zbmath_violin}.
Each execution of the algorithm requires approximately $3$ hours of computations. The algorithm seems to produce reliable conclusions concerning the top nodes visited all along the ending iterations. 
Nevertheless, it appears to be necessary to extend our investigations in order to obtain a scalable method for handling larger graphs.

%
%

\section{Proof of the main result (Theorem \ref{theo:conv})\label{sec:section_proof}}

We establish a differential inequality that will imply the convergence of $J_t$. The computations actually lead us to a system of two differential inequalities. We therefore  introduce another quantity $I_t$ that measures the average closeness (w.r.t. $x$) of the conditional law of $y$ given $x$ at time $t$ to $\nu$, defined as:
\begin{equation}\label{eq:It_def}
I_t:= \inte{KL(m_t(.|x)\vert \vert \nu) \mathrm{d} n_t(x)}= \inte{\left[ \inte{\log\dfrac{m_t(y|x)}{\nu(y)}} m_t(y|x) dy \right] \,. \mathrm{d}n_t(x)}\end{equation}

The next proposition links the evolution of $J_t$ (in terms of an upper bound of $\partial_t J_t$) with the spectral gap of $\mub$ over $\Gamma_G$, the diameter of the graph $\mathcal{D}_G$ and the divergence $I_t$.

\subsection{Study of $\partial_t J_t$}

\medskip
\begin{proposition}\label{propJ}
Let  $c^{\star}(U_{\nu})$ be the term defined in Equation \eqref{eq:cstar} and $C_{\Gamma_G}$ be the constant given in Proposition \ref{prop:log_sobolev}. We then have:
$$\partial_t J_t\le  \mathcal{D}_G^2 \beta_t' + 8 \mathcal{D}^2_G \beta_t^2 I_t- \frac{e^{-c^{\star}(U)\beta_t}}{C_{\Gamma_G}(1+\beta_t)}J_t.$$
\end{proposition}
\underline{\textit{Proof:}}
We compute the derivative of $J_t$ and separately  study  each of the three terms:

\begin{equation}\label{eq:jdecomp}
\dt J_t=\underbrace{\inte \dt\{ \log(n_t(x))\}\dn}_{J_{1,t}}-
\underbrace{\inte \dt\{ \log(\mub(x))\} \dn}_{J_{2,t}}+\underbrace{\inte\log\left[\dfrac{n_t(x)}{\mub}\right]\dt n_t(x) \de x}_{J_{3,t}}
\end{equation}
\noindent\textbf{Study of $J_{1,t}$:} This term is easy to deal with:

\begin{eqnarray}
J_{1,t} &= &\inte_{\Gamma_G} \dt\{ \log(n_t(x))\}n_t(x)  \nonumber\\
& = &\inte_{\Gamma_G} \dfrac{\dt\{ n_t(x)\}}{n_t(x)}n_t(x)=  \inte_{\Gamma_G}\dt\{ n_t(x)\} dx = \dt  \left\{ \inte_{\Gamma_G} n_t(x) dx \right\}=  \dt\{1\}=0.\label{eq:J1t}
\end{eqnarray}

\noindent\textbf{Study of $J_{2,t}$:}
Using the definition of $\mub$ given in Equation \eqref{eq:gibbs}, we obtain for the second term:

\begin{align*}
J_{2,t}&=-\inte_{\Gamma_G} \dt \{\log (\mub(x))\}\dn=-\inte_{\Gamma_G} \dt \{-\beta_t U_{\nu}(x) - \log Z_{\beta_t}\}\dn\\
&=\beta_t'\inte_{\Gamma_G} U_{\nu}(x)\dn + \frac{\dt\{Z_{\beta_t}\}}{Z_{\beta_t}}.
\end{align*}
According to the definition of $Z_{\beta_t}$, we have:
$$
\frac{\dt\{ Z_{\beta_t}\}}{Z_{\beta_t}} = Z_{\beta_t}^ {-1} \dt \left\{\inte_{\Gamma_G} e^{-\beta_t U_{\nu}(x)} \dx \right\} = - \beta_t'
 \inte_{\Gamma_G} U_{\nu}(x) \frac{e^{-\beta_t U_{\nu}(x)}}{Z_{\beta_t}} \dx = 
 - \dt\{\beta_t\} \inte_{\Gamma_G} U_{\nu}(x) \mu_{\beta_t}(x) \dx.
$$
Hence, we obtain
$$
J_{2,t} = \beta_t' \inte_{\Gamma_G} U_{\nu}(x) [n_t(x) - \mu_{\beta_t}(x)]  \dx 
$$
The graph has a finite diameter $\mathcal{D}_G$. We therefore have:  $\integ U_{\nu}(x)\dn\le \mathcal{D}_G^2$. The same inequality holds using the measure $\mub$ so that:
\begin{equation}\label{eq:J2t}
 |J_{2,t}| \le \mathcal{D}_G^2 \beta_t'.\end{equation}

\noindent\textbf{Study of $J_{3,t}$:}
The last term $J_{3,t}$ involves the backward Kolmogorov equation. 
 First, since $n_t$ is the marginal law of $X_t$, we have:
 $n_t(x)=\inte m_t(x,y) \de y$.

 Using the backward Kolmogorov equation for the Markov process $(X_t,Y_t)_{t \geq 0}$ and the Fubini theorem, we have, for any smooth enough function $f_t: x \in \Gamma_G \longmapsto \mathbb{R}$ :

\begin{align*}
\integ \fot \dt\{n_t(x)\} \de x&=\integ \fot \dt\left\{ \intev m_t(x,y) \de y \right\}\de x =\integ \intev \fot   \dt\{m_t(x,y)\} \de x \de y\\
&=\integ  \intev \mathcal{L}_t (f_t)(x) m_t(x,y)\dx  \dy
=\integ \intev [\lu+\ld](f_t)(x)  m_t(x,y) \dxy,
\end{align*}
where $\lu$ and $\ld$ are defined in Equations \eqref{eq:deflu} and \eqref{eq:defld}.
Since the function $f_t$ is independent of $y$, we have $\lu(f_t)=0$.  For the part corresponding to $\ld$, we have:
 
\begin{eqnarray*}
\lefteqn{
 \integ \intev \ld(f_t)(x) m_t(x,y) \dxy}\\
 & = & 
 \integ \intev  \left[\dfrac{1}{2}\Delta_x \fot  -\beta_t \nabla_x U_y (x)  \nabla_x f_t(x)\right] m_t(x,y) \dxy \\
 &=&\integ  \dfrac{1}{2}\Delta_x \fot n_t(x) \dx  
 -\beta_t\integ\intev \nabla_x U_y(x) \nabla_x \fot n_t(x) m_t(y|x) \dxy.
\end{eqnarray*}
because $n_t$ is the marginal distribution of $X_t$ and 
$m_t(y \vert x) \times n_t(x) = m_t(x,y).$   

\medskip

Thus, for any smooth enough function $\fot$, using the operator introduced in \eqref{eq:ltd_def}, we have:
$$\integ \fot \dt n_t(x) \de x=\integ \ltd (f_t)(x) n_t(x) \dx.$$
Replacing $\fot$ by $\lot$, we obtain:
$$J_{3,t} = \integ\log\dfrac{n_t(x)}{\mub}\dt n_t(x)=\integ \ltd \left[\log \dfrac{n_t(x)}{\mub (x)}\right] n_t(\de x).$$
Since $\mub $ is the invariant distribution of $\lcd$ (see Equation \eqref{eq:lcd_def}), it is natural to insert $\lcd$: 
\begin{align}
J_{3,t}& =
\integ \ltd \left[\log \dfrac{n_t(x)}{\mub (x)}\right] n_t(\de x)\nonumber\\
& =\integ \lcd  \left[ \log \dfrac{n_t(x)}{\mub (x)}\right] n_t(\de x)-\underbrace{\integ \left( \lcd -\ltd \right)  \left[\log \dfrac{n_t(x)}{\mub (x)}\right] n_t(\de x)}_{:=\kappa_t}.\label{eq:decomp} \end{align}

\medskip

Since
$\lcd$ is a diffusion operator and $\mub$ is its invariant measure, it is well known (see, \textit{e.g.}, \cite{BGL}) that the action of $\lcd$ on the entropy is closely linked to the Dirichlet form in the following way:

\begin{equation}\label{eq:dirichlet_form}
\integ \lcd  \left[ \log \dfrac{n_t(x)}{\mub (x)}\right] n_t(\dx)=-2\integ \left(\nabla_x \left\{\sqrt{\dfrac{n_t(x)}{\mub (x)}}\right\}\right)^2 \mub(\dx),\end{equation}
and therefore translates a mean reversion towards $\mub$ in the first term of \eqref{eq:decomp}.
\\

We now study the size of the difference between $\lcd$ and $\ltd$ and introduce the ``approximation" term of $\nabla_x\um(x)$ at time $t$:
$$\rt(x)=\intev \nabla_x (U_y)(x) m_t(y|x) \dy.$$
The relationship $\nabla_x \log(f) = 2 \frac{\nabla_x \sqrt{f}}{\sqrt{f}}$, the Cauchy-Schwarz inequality and $2 ab \leq a^2+b^2$ yield:

\begin{eqnarray*}
\left|\kappa_t\right| &= &\beta_t \left| \integ \left( \rt(x)-\nabla_x\um(x)\right)  \nabla_x\left\{ \log\dfrac{n_t(x)}{\mub (x)} \right\} n_t(\de x)\right|\\
&=&2\beta_t  \left| \integ \left( \rt(x)-\nabla_x\um(x)\right) \nabla_x\left\{ \sqrt{\dfrac{n_t(x)}{\mub (x)}} \right\} \sqrt{\dfrac{\mub(x)}{n_t(x)}}  n_t(\de x)\right| \\
&\le& 2\beta_t \sqrt{\integ \left( \rt(x)-\nabla_x\um(x)\right)^2  n_t(\de x)}\cdot \sqrt{\integ\nabla_x \left(\sqrt{\dfrac{n_t(x)}{\mub(x)}}\right)^2\mub(x)\de x}\\
&\le&\beta_t^2\inte \left( \rt(x)-\nabla_x\um(x)\right)^2 n_t(x)\de x +\integ \left(\nabla_x \left\{\sqrt{\dfrac{n_t(x)}{\mub (x)}}\right\}\right)^2 \mub(\dx)  \end{eqnarray*}

%
%
If $d_{TV}$ denotes the total variation distance, the first term of the right hand side leads to:
\begin{align*}
\left| \rt(x)-\nabla_x\um(x)\right|&=\left| \intev \nabla_x d^2(x,y)[m_t(y|x)-\nu(y)] \de y\right| \\
&\le \|\nabla_xd^2(x,y)\|_{\infty}\left|\intev [m_t(y|x)-\nu(y)]\de y\right|\\
& \le 2 \|\nabla_xd^2(x,y)\|_{\infty}d_{TV}(m_t(.\vert x),\nu)\\
&\le \sqrt{2}  \|\nabla_x d^2(x,y)\|_{\infty} \sqrt{ \intev\log \dfrac{m_t(y|x)}{\nu(y)} m_t(y|x) \de y}
\end{align*}
where the last line comes from the Csisz\'ar-Kullback inequality. Since  $d^2(\cdot,y)$ is differentiable a.e. and  its derivative is bounded for all $y\in V$ by $2 \mathcal{D}_G$, we can use $I_t$ defined in Equation \eqref{eq:It_def} to obtain:
 $$\integ \left(\rt-\rc\right)^2 (x) n_t(x) \dx \le 8 \mathcal{D}^2_G I_t $$
 Consequently, we obtain:
\begin{equation}\label{eq:boundkappa}
|\kappa_t| \leq 8 \mathcal{D}^2_G  \beta_t^2 I_t + \integ \left(\nabla_x \left\{\sqrt{\dfrac{n_t(x)}{\mub (x)}}\right\}\right)^2 \mub(\dx)
\end{equation}
 Taking  the inequalities \eqref{eq:dirichlet_form} and \eqref{eq:boundkappa}, we now obtain in \eqref{eq:decomp}:
 $$
 J_{3,t} \leq 8 \mathcal{D}^2_G \beta_t^2 I_t   - \integ \left(\nabla_x \left\{\sqrt{\dfrac{n_t(x)}{\mub (x)}}\right\}\right)^2 \mub(\dx)
 $$
 \medskip
 We denote $f_t=\sqrt{\frac{n_t}{\mub}}$. Since  $\|f\|_{2,\mub}^2=1$, one can easily see that $\enbt(f_t^2) = J_t$. Now, the Logarithmic Sobolev inequality on the (quantum) graph $\Gamma_G$ for the measure $\mub$ stated in Proposition \ref{prop:log_sobolev} shows that:
 $$\integ \left(\nabla_x \left\{\sqrt{\dfrac{n_t(x)}{\mub (x)}}\right\}\right)^2 \mub(\dx)\ge \frac{e^{-c^{\star}(U_{\nu})\beta_t}}{C_{\Gamma_G} (1+\beta_t)} J_t,$$
 where $c^{\star}(U_{\nu})$ is defined in Equation \eqref{eq:cstar} and is related to the maximal depth of a well containing a local but not global minimum.  We thus  obtain:

$$ J_{3,t}=\inte \ltd \left[\log \dfrac{n_t(x)}{\mub (x)}\right] n_t(\de x)\le 8 \mathcal{D}^2_G \beta_t^2 I_t- \frac{e^{-c^{\star}(U_{\nu})\beta_t}}{C_{\Gamma_G} (1+\beta_t)} J_t. $$
The proof is concluded by regrouping the three terms. \hfill $\square$

 

\subsection{Study of $\partial_t I_t$}

\begin{proposition}\label{propI}
Assume that $\mathcal{D}_G \geq 1$ and $\beta_t$ is an increasing inverse temperature with $\beta_t \geq 1$, then:
\begin{equation} \label{iti}
\dt I_t\le-\alpha_tI_t-\dt J_t + \mathcal{D}_G^2 [\beta_t' + 6 \beta_t^2].
\end{equation}
\end{proposition}

\underline{\textit{Proof:}}
 We compute the derivative of $I_t$. Observing that $m_t(x,y)=n_t(x)m_t(y\vert x)$, we have:
 \begin{align*}
\lefteqn{ \partial_t I_t = \partial_t \left\{\integ n_t(x) \left( \intev \log \frac{m_t(y|x)}{\nu(y)} m_t(y\vert x)\dy \right) \dx \right\}} \\
 & = \partial_t \left\{\integ  \intev \log \frac{m_t(y|x)}{\nu(y)} m_t(x,y)\dxy \right\} \\
 & = \underbrace{\integ \intev \partial_t\{\log m_t(y \vert x)\} m_t(x,y)\de x\de y}_{:=I_{1,t}}  - \underbrace{\integ \intev \partial_t\{\log \nu (y) \} m_t(x,y)\de x\de y}_{:=I_{2,t}} + \underbrace{\inte \log \dfrac{m_t(y|x)}{\nu(y)}\partial_t m_t (x,y)\de x\de y}_{:=I_{3,t}}
 \end{align*}
 The computation of the first two terms is straightforward. For the first one we have: 
 \begin{align*}
 I_{1,t}&=\integ \intev \partial_t\{ \log m_t(y|x) \}m_t(x,y)
 =\integ \intev\dfrac{\partial_t m_t(y|x)}{m_t(y|x)}m_t(x,y)\de x \de y\\
 &=\integ \intev \partial_t m_t(y|x) n_t(x)\de x\de y
 =\integ n_t(x) \left( \intev  \partial_t m_t(y|x) \de y \right) \mbox{ }\de x\\
 &=\integ n_t(x)  \partial_t \left\{\underbrace{ \intev m_t(y|x) \de y}_{:=1} \right\} \mbox{ }  \de x =0
 \end{align*}
The computation of $I_{2,t}$ is easy since $\nu$ does not depend on $t$, implying that $I_{2,t}=0$.
 
 \medskip
 \noindent
 For the third term, we use the backward Kolmogorov equation and obtain:
\begin{align*}
I_{3,t}&=\integ \intev \log\dfrac{m_t(y|x)}{\nu(y)}\partial_t m_t(x,y)=\integ \intev\mathcal{L}_t\left(\log \dfrac{m_t(y|x)}{\nu(y)} \right) m_t(x,y) \de x \mbox{ }\de y\\
&=\underbrace{\integ \intev\lu \left( \log \dfrac{m_t(y|x)}{\nu(y)}\right) m_t(x,y) \de x \mbox{ }\de y}_{:=I^1_{3,t}}+ \underbrace{\integ \intev \ld \left( \log\dfrac{m_t(y|x)}{\nu(y)}\right) m_t(x,y) \de x \mbox{ }\de y}_{:=I^2_{3,t}}
\end{align*}
The jump part $\lu$ exhibits a mean reversion on the entropy on the conditional law: applying the Jensen inequality for the logarithmic  function and the measure $\nu$, we obtain:

\begin{align}
I^1_{3,t}&= \at \integ \intev  \left[ \intev  \left(\log \dfrac{m_t(y'|x)}{\nu(y')} -\log \dfrac{m_t(y|x)}{\nu(y)} \right) \nu(y')\dy'\right] m_t(x,y)\de x\de y\nonumber\\
&=\at \integ \intev  \left[ \intev  \left( \log \dfrac{m_t(y'|x)}{\nu(y')} \right) \nu(y')\dy' \right] \mbox{ } m_t(x,y) \de x\de y-\at I_t\nonumber\\
&\le \at\integ \intev   \log \left( \intev \dfrac{m_t(y'|x)}{\nu(y')}\nu(y')\de y' \right) \mbox{ }  m_t(x,y) \de x\de y-\at I_t\nonumber\\
&\le \at\integ \intev  \log \left(\intev m_t(y'|x)\de y' \right) \mbox{ }  m_t(x,y) \de x\de y-\at I_t\nonumber\\
&\le\at\integ \intev   \log 1 \cdot m_t(x,y)\de x\de y -\at I_t\le -\at I_t\label{i31}
\end{align}
If we  consider the action of $\ld$ on the entropy of the conditional law, using  $m_t(x,y)=m_t(y\vert x) n_t(x)$ yields:
 \begin{align}
 I^2_{3,t}&=\integ \intev \ld \left( \log \dfrac{m_t(y|x)}{\nu(y)}\right) m_t(x,y) \de x \mbox{ }\de y\nonumber\\
 &=\integ \intev \ld \left( \log\dfrac{m_t(x,y)}{\nu(y)\cdot n_t(x)}\right) m_t(x,y) \de x \mbox{ }\de y\nonumber\\
 &=\integ \intev \ld \log\left(  m_t(x,y) \right) m_t(x,y) \dxy -
 \integ \intev \ld \log\left(  n_t(x,y) \right) m_t(x,y) \dxy,  \label{eq:i32dec}
 \end{align}
because $\ld(\nu)(y)=0$ ($\ld$ only involves the $x$ component). 
We now study the first term of \eqref{eq:i32dec}:

\begin{align*}
\integ \intev \ld\left(\log m_t(x,y)\right) \dm&= \dfrac{1}{2} \integ \intev \Delta_x \left[\log m_t(x,y)\right] \dm \\
&-\beta_t\integ \intev <\nabla_x U_y,\nabla_x \log m_t(x,y)>\dm \\
&=\dfrac{1}{2} \integ \intev \partial_x^2 \left\{\log m_t(x,y)\right\} \dm \\&- \beta_t \inte \partial_x U_y \partial_x \log m_t(x,y)\dm \\
 \end{align*}
The first term of the right-hand side deserves special attention:
 \begin{align*}
 \integ \intev \partial_x^2 \left\{\log m_t(x,y)\right\} \dm &= \integ \intev \partial_x \left\{\frac{\partial_x \{m_t(x,y)\}}{m_t(x,y)}\right\} \dm\\
 & = \integ \intev \frac{\partial^2_x \{m_t(x,y)\}}{m_t(x,y)} \dm  - \integ \intev  \frac{[\partial_x \{m_t(x,y)\}]^2}{m_t(x,y)} \dxy\\
 & =  \integ \intev \partial^2_x \{m_t(x,y)\}  \dxy  - 4 \integ \intev \left( \partial_x \sqrt{m_t(x,y)}\right)^2 \dxy\\
 & = \sum_{e \in E} \intev \left[ \partial_x m_t (e(L_e),y) \right]  - \partial_x m_t (e(0),y) \de y \\
 &  - 4 \integ \intev \left( \partial_x \sqrt{m_t(x,y)}\right)^2 \dxy,
  \end{align*}
  where $\{e(s),0 \leq s \leq L_e\}$ is the parametrization of the edge $e \in E$ introduced in Section \ref{sec:quantum_diffusion}. The gluing conditions \eqref{eq:glue} yield:
  $$
    \sum_{e \in E} \intev \partial_x m_t (e(L_e),y) - \partial_x m_t (e(0),y) = \intev \sum_{v \in V} \left[  \sum_{e \sim v} - d_e m_t(v,y) \right]= 0.
  $$
  As a consequence, we obtain:
  \begin{align}
  \integ \intev \ld\left(\log m_t(x,y)\right) \dm& = -2\integ \intev \left(\partial_x\sqrt{m_t(x,y)}\right)^2 \dxy \nonumber\\
  &- \beta_t \integ \intev \partial_x U_y \partial_x \log m_t(x,y)\dm \label{eq:boundL2}
  \end{align}
 The Cauchy-Schwarz inequality  applied on the second term leads to:
\begin{align*}
\left| \beta_t \integ\intev \partial_x U_y \partial_x \log m_t(x,y)\dm  \right|&\le \beta_t \|\partial_x d^2(.,.)\|_{2,m_t}\sqrt{\integ \intev \left(\partial_x \log m_t(x,y)\right)^2 \dm} \\
&\le \beta_t ||\partial_x d^2(x,y)||_{\infty} \sqrt{\integ \intev \left(\partial_x \log m_t(x,y)\right)^2 m_t(x,y)\de x \de y}\\
&\le 4 \beta_t \mathcal{D}_G \sqrt{\integ \intev \left(\partial_x \sqrt{m_t(x,y)}\right)^2 \de x \de y}\\
&\le 2\mathcal{D}^2_G \beta_t^2  + 2\integ \intev \left(\partial_x \sqrt{m_t(x,y)}\right)^2  \de x\de y\\
\end{align*} 
Inserting this in Equation \eqref{eq:boundL2} leads to:

\begin{equation}\label{eq:ldlog}
 \integ \intev \ld\log m_t(x,y)\dm\le  2\mathcal{D}_G^2 \beta_t^2 .\end{equation}
We study the second term of \eqref{eq:i32dec} and use the proof of Proposition $\ref{propJ}$: 
 the decomposition  \eqref{eq:jdecomp} with Equations \eqref{eq:J1t} and \eqref{eq:J2t} yield:
 $$\partial_t J_t\le \beta_t' \mathcal{D}_G^2+ \inte \ld \left( \lot \right) m_t(\de x,\de y). $$
This implies:
\begin{equation}\label{eq:interm}-\integ \intev \ld \left( \log n_t(x)\right) m_t(\dx,\dy)\le -\dt J_t + \beta_t' \mathcal{D}_G^2-\inte \ld \left( \log\mub(x)\right) m_t(\de x, \dy).\end{equation}
Using the definition of $\mub$, we obtain:
\begin{align}
-\integ \intev \ld \log \dfrac{e^{-\beta_t U_{\nu}(x)}}{Z_{\beta_t}}&=\integ \intev \left[\ld \left(  \beta_t U_{\nu}(x)\right)+\ld  \left( \log \z \right) \right]m_t(\de x,\de y)\nonumber\\
& = \beta_t  \integ \intev \ld \left( U_{\nu}(x)\right) m_t(\de x,\de y)\nonumber\\
&=\beta_t \integ \intev  \dfrac{1}{2}\Delta_x \um (x)-\beta_t\nabla_x U_y(x) \nabla_x U_{\nu}(x) m_t(\de x,\de y)\nonumber\\
&\le \frac{\beta_t}{2} + 4 \beta_t^2 \mathcal{D}_G^2\nonumber
\end{align}
This inequality used in \eqref{eq:interm} yields:
\begin{align}
-\inte \ld\left( \log n_t(x) \right) m_t(x,y)&=-\dt J_t + \beta_t' \mathcal{D}_G^2 +  \frac{\beta_t}{2} + 4 \beta_t^2 \mathcal{D}_G^2\label{eq:boundpartialJ}
\end{align}
We now use \eqref{eq:ldlog} and \eqref{eq:boundpartialJ} in \eqref{eq:i32dec} and our assumptions $\beta_t \geq 1$ and $\mathcal{D}_G \geq 1$ to obtain: 
\begin{equation}\label{i32}
I^2_{3,t}\le  -\dt J_t + \mathcal{D}_G^2 [\beta_t' + 6 \beta_t^2]
\end{equation}
Combining $\eqref{i31}$ and $\eqref{i32}$ leads to the desired inequality given by Equation \eqref{iti}.
\hfill $\square$

\subsection{Convergence of the entropy} 
The use of Propositions \ref{propJ} and \ref{propI} makes it possible to obtain a system of coupled differential inequalities.\\

\noindent
\underline{\textit{Proof of Theorem \ref{theo:conv}}:}
 If we denote $a = \mathcal{D}_G^2$, we can write:
$$\begin{cases}
J'_t\le - \frac{e^{-c^{\star}(U_{\nu})\beta_t}}{C_{\Gamma_G} (1+\beta_t)}J_t + a \beta_t' + 8 a \beta_t^2 I_t  \\
I'_t\le-\alpha_tI_t-J'_t +a (\beta'_t +6 \beta_t^2)\end{cases}
$$
We introduce an auxiliary function $K_t=J_t+k_t I_t$, where $k_t$ is a smooth positive decreasing function for $t$ large enough, so that $\displaystyle\lim_{t\to \infty}k_t=0$. We use the system above to deduce that:
\begin{align*}
 K'_t&=J'_t+ k_t' I_t+ k_t I'_t\\
&\le  J'_t+ k_t I'_t\\
&\le  J'_t + k_t(-\alpha_t I_t- J'_t)  + a k_t (\beta'_t +6 \beta_t^2)\\
&\le (1-k_t) J'_t- k_t\alpha_t I_t + a k_t (\beta'_t +6 \beta_t^2).
\end{align*}
In the first inequality, we use the fact that  $k'_t\leq 0$ and $I_t$ is positive. The second inequality is given by $ \eqref{iti}$. The upper bound of Proposition \ref{propJ} now leads to:
$$
K'_t \leq - \epsilon_t (1-k_t) J_t  - k_t \alpha_t I_t + 8 a (1-k_t) \beta_t^2 I_t + a \beta_t' +  6 a k_t  \beta_t^2,
$$
where we denoted $\epsilon_t = \frac{e^{-c^{\star}(U_{\nu})\beta_t}}{C_{\Gamma_G} (1+\beta_t)}$. 
We choose the function $k_t$ to obtain a mean reversion on $K_t$:
\begin{equation}
k_t := \frac{8 a \beta_t^2}{\alpha_t + 8 a \beta_t^2 - \epsilon_t/2}.
\end{equation}
Note that this function is decreasing for sufficiently large $t$ as soon as $\beta_t = o(\alpha_t)$, which is the case according to the choices described in Theorem \ref{theo:conv}. Moreover, a straightforward consequence is $\lim_{t \longrightarrow + \infty} k_t = 0$.
This ensures that a positive $T_0$ exists such that:
$$
\forall t \geq T_0 \qquad 0 \leq k_t \leq \frac{1}{2}.
$$
Consequently, we deduce that:
\begin{align*}
\forall t \geq T_0 \qquad 
K_t'& \leq - \epsilon_t (1-k_t) J_t -  k_t  \frac{\epsilon_t}{2} I_t+ a \beta_t' + 6 a \beta_t^2 k_t\\
& \leq  - \frac{\epsilon_t}{2} J_t -  k_t  \frac{\epsilon_t}{2} I_t+ a \beta_t' + 6 a \beta_t^2 k_t \leq - \frac{\epsilon_t}{2} K_t  + \underbrace{a \beta_t' + 6 a \beta_t^2 k_t}_{:=\eta_t}\\
\end{align*}

The next bound is an easy consequence of the Gronwall Lemma:
$$
\forall t \geq T_0 \qquad 
K_t \leq \displaystyle K_{T_0} e^{-\displaystyle\int_{T_0}^t \frac{\epsilon_s}{2} \de s} + \displaystyle\int_{T_0}^T \eta_s e^{-\displaystyle\int_{s}^t \frac{\epsilon_u}{2} \de u} \de s,
$$
which in turn implies that $K_t \longrightarrow 0$ as $t \longrightarrow + \infty$ as soon as
$\eta_t = o_{t \sim + \infty} (\epsilon_t)$ with $\int_{0}^{\infty} \epsilon_u du = + \infty$.

\medskip

We are now looking for a suitable choice for $\alpha_t$ and $\beta_t$. Let us assume that $\alpha_t \sim \lambda t^{\gamma}$ for any $\gamma>0$. This choice leads to:
$$k_t \sim  \frac{8 a b^2  \log(t+1)^2}{\lambda t^{\gamma}},$$ so that:
$$
\eta_t \sim \frac{ab}{t} + \frac{48 a^2 b^4 \log(t+1)^4}{\lambda t^{\gamma} } = O(t^{-(1 \wedge \gamma)} \log(t)^4).
$$
At the same time, we can check that 
$\epsilon_t \sim \frac{t^{-c^{\star}(U_{\nu}) b }}{b C_{\Gamma_G} \log t}$.
Now, our conditions on $\eta_t$ and $\epsilon_t$ imply that:
$$
\int_{0}^{\infty} \epsilon_u \de u= +\infty \qquad \Longleftrightarrow \qquad b \, c^\star(U_\nu) \leq 1.
$$
At the same time: 
$$
 \eta_t = o(\epsilon_t) \qquad \Longleftrightarrow \qquad  1 \wedge \gamma > b  \, c^\star(U_\nu).
$$
The optimal calibration of our parameters $\gamma$ and $b$ (minimal value of $\gamma$, maximal value of $b$) induces the choice $\gamma=1$ and $b < c^\star(U_\nu)^{-1}$.
Up to the choices $\alpha_t \sim \lambda t$ and $\beta_t \sim b \log(t)$, we deduce that 
 $K_t \longrightarrow 0$  when $t$ goes to infinity. 
  Since $I_t$, $k_t$ and $J_t$ are positive, this also implies that $J_t \longrightarrow 0$  when $t$ goes to infinity. 
\hfill $\square$

\section{Functional inequalities}\label{sec:functional}
This section is devoted to the proof of the Log-Sobolev inequality for the measure $(\mu_{\beta})_{\beta >0}$. We dealt with the Dirichlet form using this inequality (see Equation \eqref{eq:dirichlet_form})    in Proposition \ref{propJ}.
In particular, we need to obtain an accurate estimate when $\beta \longrightarrow +\infty$.
For this purpose, we introduce the generic notation for  Dirichlet forms (see, \textit{i.e.}, \cite{BGL} for a more in-depth description):
$$\forall \beta \in \mathbb{R}_+^* \quad \forall f \in W^{1,2}(\mubb) \qquad \eb(f,f)= \|\grad f\|^2_{\mubb} = \integ | \grad f(x)|^2 \de \mubb(x).$$
If we denote $<f>_{\mubb} = \mubb(f) = \integ f \de \mubb$, we are interested in showing the Poincar\'e inequality:
\begin{equation}\label{eq:pi}
\|f-<f>_{\mubb}\|_{2,\mubb}^2 = \integ ( f-<f>_{\mubb})^2 \de \mubb \leq \lambda(\beta) \eb(f,f) 
\end{equation}
and the Log-Sobolev Inequality (referred to as LSI below):
\begin{equation}\label{eq:lsi}
\int_{\Gamma_G} f^2 \log \left( \frac{f}{\|f\|_{2,\mubb}}\right)^2 \de \mubb \leq C(\beta) \eb(f,f)
\end{equation}

The specific feature of this functional inequality deals with the quantum graph settings and deserves careful adaptation of the pioneering work of \cite{stroock_sa1}. This technical section is split into two parts. The first one establishes a preliminary estimate when $\beta = 0$ (\textit{i.e.}, when dealing with the uniform measure on $\Gamma_G$). The second one then uses this estimate to derive the asymptotic behavior of the LSI when $\beta \longrightarrow +\infty$.

\subsection{Preliminary control for $\mu_{0}$ on $\Gamma_G$}

We consider $\mu_0$ the normalized Lebesgue measure and use the standard notation for any measure $\mu$ on $\Gamma_G$:
$$
\|f\|_{W^{1,p}(\mu)} := \|f\|_{p,\mu}+\|\nabla f\|_{p,\mu}
$$

Let us establish the next elementary result: 
\begin{lemma}\label{lemma:sobolev_mu0} The ordinary Sobolev inequality holds on $\Gamma_G$, i.e., for all measurable functions $f$ we have:

\begin{equation}
\|f- <f>_0 \|_{p,\mu_0}^2\le e^{2/e} L^2 \left[ \frac{\mathcal{D}_G L}{2}+1\right]   \mathcal{E}_{0}(f,f),
\end{equation}
where $\mathcal{D}_G$ is the diameter of the graph and $L$ its perimeter defined by:
$$
L = \sum_{e \in E} L_e.
$$
\end{lemma}

\underline{\textit{Proof:}}
First, let us remind the reader that for a given interval $I$ in dimension 1 equipped with the Lebesgue measure $\lambda_I$, the Sobolev space $W^{1,p}(\lambda_I)$ is continuously embedded  in $L^{\infty}(I)$ (compact injection when $I$ is bounded). In particular, it can be shown (see, \textit{e.g.}, \cite{Brezis}) that while integrating w.r.t. the unnormalized Lebesgue measure:
$$ \forall p \geq 1 \qquad 
\| f\|_{L^{\infty}(I)}\le e^{1/e} \|f\|_{W^{1,p}(\lambda_I)}.  $$


We now consider $f \in W^{1,p}(\mu_0)$. Since $G$ has a finite number of edges we can write: $\Gamma _G= \bigcup_{i=1}^n e_i,$ where each $e_i$ can be seen as an interval of length $\ell_i$. We have seen that  for each edge $e_i$:
$$\| f\|_{L^{\infty}(e_i)}\le e^{1/e} \|f\|_{W^{1,p}(\lambda_{e_i})}.$$
We can use a union bound since $\Gamma_G$ represents the union of edges, and deduce that:
$$\| f\|_{L^{\infty}(\Gamma_G)}=\max_{1 \leq i \leq n}\| f\|_{L^{\infty}(e_i)} \leq e^{1/e} \max_{1 \leq i \leq n}\|f\|_{W^{1,p}(\lambda_{e_i})} \leq e^{1/e} \|f\|_{W^{1,p}(\lambda_{\Gamma_G})},$$
where the inequality above holds w.r.t. the unnormalized Lebesgue measure. If $L$ denotes the sum of the lengths of all edges in $\Gamma_G$, we obtain:
\begin{equation}\label{eq:injection}\| f\|_{L^{\infty}(\Gamma_G)}  \leq e^{1/e} L \|f\|_{W^{1,p}(\mu_0)}.\end{equation}

\medskip

Second, we establish a simple Poincar\'e inequality for $\mu_0$ on $\Gamma_G$. For any function $f \in W^{1,2}(\mu_0)$, we use the equality:
$$
\|f-<f>_0\|_{2,\mu_0}^2 = \frac{1}{2} \integ \integ [f(x)-f(y)]^2 \de \mu_0(x) \de \mu_0(y) = \frac{1}{2}
\integ \integ  \left[ \int_{\gamma_{x,y}}f'(s) ds \right] ^2 \de \mu_0(x) \de \mu_0(y).
$$
where $\gamma_{x,y}$ is the shortest path that connects $x$ to $y$, parametrized with speed $1$ and $f'(s)$ refers to the derivative of $f$ w.r.t. this parametrization at time $s$. It should be noted that such a path exists because the graph $\Gamma_G$ is connected. The Cauchy-Schwarz inequality yields:
\begin{equation}\label{eq:poincare}
\|f-<f>_0\|_{2,\mu_0}^2 \leq \frac{1}{2}
\integ \integ  \left[ \int_{\gamma_{x,y}}|\nabla f(s)|^2 ds \right] |\gamma_{x,y}| \de \mu_0(x) \de \mu_0(y) \leq \frac{\mathcal{D}_G L }{2} \|\nabla f\|_{2,\mu_0}^2.
\end{equation}

\medskip 

The Sobolev inequality is now an obvious consequence of the previous inequality: consider $f \in W^{1,p}(\mu_0)$ and note that since
$\mu_0(\Gamma_G)=1$, then \eqref{eq:injection} applied with $p=2$ leads to:
\begin{eqnarray*}
\forall q \geq 1 \qquad 
\|f-<f>_0\|_{q,\mu_0}^2  & =& \left(\integ |f-<f>_0|^q \de \mu_0\right)^{2/q}\\
& \leq& \|f-<f>_0\|_{L^{\infty}(\Gamma_G)}^2 \mu_0(\Gamma_G)\\
& \leq &e^{2/e} L^2 \|f-<f>_0\|^2_{W^{1,2}(\mu_0)}.\\
\end{eqnarray*}
We can use the Poincar\'e inequality established for $\mu_0$ on $\Gamma_G$ in Equation \eqref{eq:poincare} and obtain:
$$
\forall q \geq 1 \qquad 
\|f-<f>_0\|_{q,\mu_0}^2 \leq e^{2/e} L^2 \left[ \frac{\mathcal{D}_G L}{2}+1\right]  \|\nabla f\|_{2,\mu_0}^2 = e^{2/e} L^2 \left[ \frac{\mathcal{D}_G L}{2}+1\right]   \mathcal{E}_{0}(f,f),
$$
which concludes the proof. 
\hfill $\square$

\begin{remark}
The constants obtained in the proof of Lemma \ref{lemma:sobolev_mu0} above could certainly be improved. Nevertheless,  such an  improvement would have little importance for the final estimates obtained in Proposition \ref{prop:log_sobolev}.
\end{remark}

\subsection{Poincar\'e Inequality on  $\mu_{\beta}$}\label{sec:poincare_section}

In the following, we show a Poincar\'e Inequality for the measure $\mu_{\beta}$ for large values of $\beta$. This preliminary estimate will be useful for deriving LSI on $\mubb$. This functional inequality is strongly related to the classical minimal elevation of the energy function $U_\nu$ for joining any state $x$ to any state $y$.

We first introduce  some useful notations. For any couple of vertices $(x,y)$ of $\Gamma_G$, and for any  path  $\gamma_{x,y} $  that connects them, we define $h(\gamma_{x,y})$ as the highest value of $U_{\nu}$ on $\gamma_{x,y}$: 
$$h(\gamma_{x,y}) = \max_{s\in \gamma_{x,y}} U_{\nu}(s) .$$ 
We define $H(x,y)$ as the smallest value of $h(\gamma_{x,y})$ obtained for all possible paths from $x$ to $y$:
$$ H(x,y)=\min_{\gamma: x\to y} h(\gamma)$$ 
Now, for any pair of vertices $x$ and $y$, the notation $\gamma_{x,y}$ will be reserved for the path that attains the minimum in the definition of $H(x,y)$. Such a path exists for any $x,y$ because $\Gamma_G$ is connected and possesses a finite number of paths that connect any two given vertices. 

Finally, we introduce the quantity that will mainly determine the size of the spectral gap involved in the Poincar\'e inequality and the constant in the LSI (see the seminal works of \cite{FWFW} for a LDP probabilistic interpretation and \cite{stroock_sa} for a  functional analysis point of view):
\begin{equation}\label{eq:cstar}
c^{\star}(U_{\nu}) : = \max_{(x,y) \in \Gamma_G^2}\left[ H(x,y)-U_\nu(x)-U_\nu(y)\right]+\min_{x\in \Gamma_G}U_\nu(x)
\end{equation}

In the following, every time we write $\min U_{\nu}$ we refer to $\displaystyle{\min_{x\in \Gamma_G}}U_\nu(x)$.

\begin{theorem}[Poincar\'e inequality for $\mubb$]\label{prop:poincare_beta}
For all measurable functions  $g$ defined on $\Gamma_G $:
$$ Var_{\mubb}(g)\le \frac{|E| \max_{e \in E} L_e  e^{2 \mathcal{D}_G}}{2} e^{\beta c^{\star}(U_{\nu}) } \mathcal{E}_{\mubb}(g,g).$$

\end{theorem}

\underline{\textit{Proof:}}
Since the graph $\Gamma_G$ is connected, for any two points $x$ and $y$, we can find a minimal path $\gamma_{x,y}$ that links them and that minimizes $H$.
 We denote $x_0=x,x_1 ,\cdots x_{l},x_{l+1}=y$ where the sequence $x_1,\cdots,x_l$ refers to the nodes included in the path $\gamma_{x,y}$.
\begin{align*}
Var_{\mubb}(g)&=\dfrac{1}{2} \integ \integ (g(y)-g(x))^2\de \mubb(x) \de \mubb(y)\\
&=\dfrac{1}{2} \integ \integ \left(\sum_{i=0}^{l}(g(x_{i+1})-g(x_i)\right)^2\de \mubb(x) \de \mubb(y)\\
&=\dfrac{1}{2}\integ \integ \left(\sum_{i=0}^{l}\displaystyle\int_{x_i}^{x_{i+1}} g'(s)\de s\right)^2\de \mubb(x) \de \mubb(y)\\
&=\dfrac{1}{2}\integ\integ\left( \intega g'(s)\de s\right)^2\de \mubb(x) \de \mubb(y)\\
&\le \dfrac{1}{2} \integ\integ \intega | \nabla g(s)|^2 \de s |\gamma_{x,y}|\de \mubb(x) \de \mubb(y),\\
\end{align*} 
The last line is implied by the Cauchy-Schwarz inequality.

\medskip
We can organize the terms involved in the above upper bound in the following way:
for any  edge of the graph $e \in E$, we denote $V_e$  the set of points $(x,y)$ such that $ \gamma_{x,y} \cap e \neq \emptyset$.
We therefore have:
\begin{align*}
Var_{\mubb}(g)&\le \dfrac{1}{2} \sum_{e\in E}  \intee | \nabla g(s)|^2 \de s \inteve |\gamma_{x,y}|\de \mubb(x) \de \mubb(y)\\
&\le  \dfrac{1}{2} \sum_{e\in E}  \intee | \nabla g(s)|^2 \mubb(s)  \inteve |\gamma_{x,y}|\dfrac{1}{\mubb(s)}\de \mubb(x) \de \mubb(y)\de s\\
&\le  \dfrac{1}{2} \sum_{e\in E} \left[ \intee  | \nabla g(s)|^2 \mubb(s) \de s\right] \left[\mbox{  } \sup_{\tilde{s}\in e}\dfrac{1}{\mubb(\tilde{s})} \inteve |\gamma_{x,y}|\de \mubb(x) \de \mubb(y)\right]\\
&\le  \dfrac{1}{2}   \integ   | \nabla g(s)|^2 \mubb(s) \de s\mbox{  } \sup_{\tilde{s}\in \Gamma_G}\dfrac{1}{\mubb(\tilde{s})} \inteves |\gamma_{x,y}|\de \mubb(x) \de \mubb(y)\\
\end{align*}
In this case,  $e_{\tilde{s}}$ denotes the edge of the graph that contains the point $\tilde{s}$ and the set $V_{e_{\tilde{s}}}$ is still the set of couples $(x,y)$ defined  above, associated with each edge $e_{\tilde{s}}$.
%
%
%
%
%
%
%
%
%
We introduce the quantity $\mathcal{A}$ defined as:
$$ \mathcal{A}=\sup_{\tilde{s}\in \Gamma_G}\dfrac{1}{\mubb(\tilde{s})} \inteves |\gamma_{x,y}|\de \mubb(x) \de \mubb(y).$$
Using this notation, we have obtained that for all functions $g$, we have the Poincar\'e inequality:
\begin{equation}\label{eq:poinc}
 Var_{\mubb}(g)\le \frac{\mathcal{A}}{2} \mathcal{E}_{\mubb}(g,g).\end{equation}

\medskip
\noindent
All that remains to be done is to obtain an upper bound of $\mathcal{A}$.
\begin{align*}
\mathcal{A}&=\sup_{\tilde{s}\in \Gamma_G}\dfrac{1}{\mubb(s)} \inteves |\gamma_{x,y}|\de \mubb(x) \de \mubb(y)\\
&=\sup_{\tilde{s}\in \Gamma_G} \zb e^{\beta U_{\nu}(\tilde{s})}\inteves \dfrac{e^{-\beta  (U_{\nu}(x)+U_{\nu}(y))}}{\zb^2}|\gamma_{x,y}| \de x\de y\\
&=\dfrac{1}{\zb}\sup_{\tilde{s}\in \Gamma_G}  \inteves e^{\beta (U_{\nu}(\tilde{s})-U_{\nu}(x)-U_{\nu}(y))}|\gamma_{x,y}| \de x\de y\\
\end{align*}
Since $\gamma_{x,y}$ is the minimal path for $H(x,y)$, we have  $H(x,y)=\displaystyle\max_{s\in\gamma_{x,y}} U_{\nu}(s)$.
Therefore:
\begin{align*}
\mathcal{A}&\le \dfrac{1}{\zb}\sup_{\tilde{s}\in \Gamma_G}  \inteves e^{\beta ((H(x,y)-U_{\nu}(x)-U_{\nu}(y))}|\gamma_{x,y}| \de x\de y\\
&\le \dfrac{1}{\zb}\sup_{\tilde{s}\in \Gamma_G}  \inteves e^{\beta (c^\star(U_{\nu})-\min(U_{\nu}))}|\gamma_{x,y}| \de x\de y\\
&\le \dfrac{e^{\beta ( c^\star(U_{\nu})-\min(U_{\nu}))}}{\zb}\sup_{\tilde{s}\in \Gamma_G}  \inteves |\gamma_{x,y}| \de x\de y\\
&\le  \dfrac{e^{\beta ( c^\star(U_{\nu})-\min(U_{\nu}))}}{\zb} |E| \max_{e \in E} L_e.
\end{align*} 
Using the definition of $\zb$, we have:
\begin{equation}\label{inegA}
\mathcal{A} \le  \dfrac{ e^{\beta c^\star(U_{\nu})}|E| \max_{e \in E} L_e}{\integ e^{-\beta(U_{\nu}(x)-\min(U_{\nu}))} \dx}.
\end{equation}
If $x^\star \in M_ {\nu}$ is a Fr\'echet mean that minimizes $U_{\nu}$, we designate $\bule$, as the ball of center $x^\star$ and radius $\dfrac{1}{\beta}$, for the geodesic distance $d$ on the graph $\Gamma_G$.
It is easy to check that:
$$
|U_{\nu}(x)-U_{\nu}(x^\star)| = \left|E_{Y\sim \nu}\left[ d^2(x,Y) - d^2(x^\star,Y) \right]\right| \leq E_{Y\sim \nu}\left| d^2(x,Y) - d^2(x^\star,Y) \right| \leq 2 \mathcal{D}_G \times d(x,x^\star).
$$
We can then deduce a lower bound on the denominator involved in \eqref{inegA}:
\begin{align*}
\integ e^{-\beta(U_\nu(x)-\min(U_\nu))}\de x&\ge \inte_{\bule} e^{-\beta(U_\nu(x)-\min(U_\nu))}\de x\\
&\ge \inte_{\bule} e^{-2 \mathcal{D}_G \beta   d(x,x^\star)} \de x\\
&~\ge e^{-2 \mathcal{D}_G} \inte_{\bule}  \de x .
\end{align*}
The Lebesgue measure of $\bule$ may be lower bounded by $\beta^{-1}$ since there is, at the least, one path in $\Gamma_G$ passing by the point $x^\star$. Inserting this inequality in  \eqref{inegA} gives:
$$ \mathcal{A} \le |E| \max_{e \in E} L_e e^{2 \mathcal{D}_G}\,  \beta e^{\beta c^\star(U_\nu)}.$$
Using this upper bound in \eqref{eq:poinc} leads to the desired Poincar\'e inequality.
\hfill $\square$

\subsection{Sobolev Inequalities  on  $\mu_{\beta}$}\label{sec:sobolev_section}

\subsubsection{Preliminary control on Dirichlet forms}

The next result will be useful to derive a LSI for $\mu_{\beta}$ from a Poincar\'e inequality on $\mu_{\beta}$ (given in Theorem \ref{prop:poincare_beta}). It generalizes the Poincar\'e inequality for $p$ norms with $p>2$ while using the Sobolev inequality given in Lemma \ref{lemma:sobolev_mu0}.

First, we introduce the maximal elevation of $U_{\nu}$ as:
\begin{equation}
\label{eq:Delta}
M=\displaystyle\sup_{x\in \Gamma_G} U_{\nu}(x) -\inf_{x\in  \Gamma_G}U_{\nu}(x).\end{equation}
Note that in our case, $M$ may be upper bounded by $\mathcal{D}_G^2$.
\begin{proposition}\label{prop:pnorm}
For any $p>2$ and all measurable functions $f$, we have:
 \begin{equation}\label{eq:prelim}
\forall \beta \geq 0 \qquad  
 \|f- <f>_{\beta} \|_{p,\mubb}^2\le 4 L^2 \left( \frac{\mathcal{D}_G L}{2}+1 \right)  e^{2/e+M \beta} \mathcal{E}_{\beta} (f,f) .
\end{equation}
\end{proposition}

\underline{\textit{Proof:}}
 The Jensen inequality applied to the convex function $x\longmapsto x^p$ yields:
\begin{align*}
\|f-<f>_{\beta}\|_{p,\mubb}^p&=\integ |f-<f>|_{\beta}^p \de \mubb \\
&=  \integ \left|f(x)-\integ f(y) \de \mubb(y)\right| ^p \de \mubb(x) \\
&=  \integ \left| \integ f(x)-f(y) \de \mubb(y)\right| ^p \de \mubb(x) \\
&\le   \integ  \integ\left|  f(x)-f(y)\right| ^p \de \mubb(y) \de \mubb(x) \\
&\le   \integ\left(  |f(x)-<f>_0| + |f(y)-<f>_0| \right) ^p \de \mubb(y) \de \mubb(x)\\
\end{align*}
Again, the Jensen inequality $|a+b|^p \leq 2^{p-1} [|a|^p+|b|^p]$ implies that:
\begin{align*}
\|f-<f>_{\beta}\|_{p,\beta}^p&\le   2^{p-1} \integ  |f(x)-<f>_0|^p + |f(y)-<f>_0|^p  \de \mubb(y) \de \mubb(x) \\
&\le 2^p \integ  |f(x)-<f>_0|^p \de \mubb(x) = 2^p \|f-<f>_0\|_{p,\mubb}^p \\
\end{align*}
We conclude that:
 \begin{equation}\label{eq:shift_beta} \|f-<f>_{\beta}\|_{p,\mubb}^2\le 4 \|f-<f>_{0}\|_{p,\mubb}^2.\end{equation}

Using the fact that:
$$
\mu_{\beta}(x) = \frac{e^{-\beta \min U_{\nu}}}{Z_\beta} \leq \frac{e^{-\beta \min U_{\nu}}}{Z_\beta} \times Z_{0} \mu_{0}(x),
$$
we also have:
$ \|f-<f>_{0}\|_{p,\mubb}^p\le \dfrac{L e^{-\beta \min U_{\nu}}}{\zb}\|f-<f>_{0}\|_{p,\mu_0}^p,$ since $Z_0$ is the perimeter $L$ of the graph $\Gamma_G$. Consequently, we obtain:
$$
 \|f-<f>_{0}\|_{p,\mubb}^2 =  \left(\|f-<f>_{0}\|_{p,\mubb}^p\right)^{2/p} \leq \left( \dfrac{L e^{-\beta \min U_{\nu}}}{\zb}\right)^{2/p}\|f-<f>_{0}\|_{p,\mu_0}^2
$$
Using inequality $\eqref{eq:shift_beta}$, the fact that $\zb \leq L e^{-\beta \min U_{\nu}}$ and the assumption $2/p<1$, we conclude that:
$$
 \|f-<f>_{\beta}\|_{p,\mubb}^2\le \dfrac{4 L e^{-\beta \min U_{\nu}}}{\zb}\|f-<f>_{0}\|_{p,\mu_0}^2.
$$
The Sobolev inequality \eqref{lemma:sobolev_mu0} implies that:
\begin{equation}\label{norme p,beta < e0}
\|f-<f>_{\beta}\|_{p,\mubb}^2\le \dfrac{4 L^3 \left( \mathcal{D}_G L/2+1\right) e^{2/e -\beta \min U_{\nu}}}{\zb}\mathcal{E}_{0}(f,f).
\end{equation}
We now  find an upper bound for the Dirichlet form $\mathcal{E}_{0}(f,f)$ that involves
$\mathcal{E}_{\beta}(f,f)$:
\begin{align}
\mathcal{E}_0(f,f)&=\integ<\grad f,\grad f> \de \mu_0\nonumber\\
&=\zb \integ <\grad f,\grad f>e^{\beta U_{\nu}(x)}  \dfrac{e^{-\beta U_{\nu}(x)}}{\zb}  \de \mu_0 \nonumber\\
&\le \displaystyle\frac{\zb \exp \left(\beta \displaystyle\sup_{x\in \Gamma_G}U_{\nu}(x)\right)}{Z_0} \integ <\grad f,\grad f>\de \mubb\nonumber\\
&\le  \displaystyle\frac{\zb \exp \left(\beta \displaystyle\sup_{x\in \Gamma_G}U_{\nu}(x)\right)}{L} \eb(f,f) \label{ineg}
\end{align}
Putting  $\eqref{norme p,beta < e0}$ and $\eqref{ineg}$ together concludes the proof. 
\hfill $\square$

\subsubsection{Log-Sobolev Inequality on $\mubb$}
For all probability measures $\mu$ and all measurable functions $f$, we denote:
$$ \en(f^2)  = \integ f^2 \log\left(\dfrac{f^2}{\|f\|^2_{2,\mu}}\right) \de \mu$$
\begin{proposition}\label{prop:log_sobolev}
The Log-Sobolev Inequality holds on $\Gamma_G$. A constant $C_{\Gamma_G} $ exists such that for all $\beta\ge 0$ and all $\mubb$-measurable functions $f$, we have:  
$$\enb(f^2)\leq  C_{\Gamma_G}[1+\beta] e^{c^\star(U_{\nu}) \beta} \eb(f,f) 
$$
\end{proposition}
\textit{\underline{Proof:}}
We consider $\beta \geq 0$ and a $\mubb$ measurable function $f$. We apply the Jensen inequality for the logarithmic function and the measure $f^2/\|f\|^2_{2,\mubb} \de \mubb$ to obtain:

\begin{align*}
\integ \left(\dfrac{f}{\|f\|_{2,\mubb}}\right)^2\log \left(\dfrac{f^2}{\|f\|^2_{2,\mubb}}\right) \de \mubb &=\dfrac{2}{p-2}\integ \left(\dfrac{f}{\|f\|_{2,\mubb}}\right)^2\log \left(\dfrac{|f|^{p-2}}{\|f\|^{p-2}_{2,\mubb}}\right) \de \mubb\\
&\le \dfrac{2}{p-2}\log\left(\integ \dfrac{|f|^{p-2}}{\|f\|_{2,,\mubb}^{p-2}}  \dfrac{f^2}{\|f\|_{2,,\mubb}^{2}}\de \mubb \right)\\
&\le \dfrac{2}{p-2}\log\left(\dfrac{\|f\|^p_{p,\mubb}}{\|f\|^p_{2,\mubb}}\right)=\dfrac{p}{p-2}\log\left(\dfrac{\|f\|^2_{p,\mubb}}{\|f\|^2_{2,\mubb}}\right).
\end{align*}
Observing that for all  $x, \delta>0$  we have $\log (x \delta) \le \delta x$, we obtain 
$\log(x) \leq x \delta + \log(\delta^{-1}).$ Therefore:

\begin{align*}
\forall \delta >0 \qquad 
\integ f^2 \log \left(\dfrac{f}{\|f\|_{2,\mubb}}\right)^2 \de \mubb &=\|f\|_{2,\mubb}^2  \integ \left(\dfrac{f}{\|f\|_{2,\mubb}}\right)^2\log \left(\dfrac{f^2}{\|f\|^2_{2,\mubb}}\right) \de \mubb \\
&\le  \dfrac{p \|f\|_{2,\mubb}^2}{p-2}\log \left(\dfrac{\|f\|^2_{p,\mubb}}{\|f\|^2_{2,\mubb}}\right)\\
&\le  \dfrac{p \|f\|_{2,\mubb}^2 }{p-2} \left[ \delta \dfrac{\|f\|^2_{p,\mubb}}{\|f\|^2_{2,\mubb}} + \log\dfrac{1}{\delta} \right].
\end{align*}
Replacing $f$ with  $f-<f>_{\beta}$ and choosing $\delta = e^{-M\beta}$ leads to:
\begin{align*}
\integ (f-<f>_{\beta})^2 &\log\left(\dfrac{f-<f>_{\beta}}{\|f-<f>_{\beta}\|_{2,\beta}}\right)^2 \de \mubb \\
&\le \dfrac{p}{p-2}\left[e^{-\beta M}\|f-<f>_{\beta}\|_{p,\mubb}^2+ M\beta \|f-<f>_{\beta}\|_{2,\mubb}^2 \right].
\end{align*}
Proposition \ref{prop:pnorm} and the Poincar\'e inequality established in Theorem \ref{prop:poincare_beta} yields: 

\begin{eqnarray}
\enb \left[ (f-<f>_{\beta})^2 \right] & = & 
\integ (f-<f>_{\beta})^2 \log\left(\dfrac{f-<f>_{\beta}}{\|f-<f>_{\beta}\|_{2,\mubb}}\right)^2 \de \mubb \nonumber\\
&\le & \dfrac{p}{p-2}\left(4 L^2 \frac{\mathcal{D}_G L +2 }{2} e^{2/e} + M \beta |E| \max_{e \in E} L_e \frac{e^{2 \mathcal{D}_G}}{2} \beta e^{c^\star(U_{\nu}) \beta} \right)\eb(f,f)\label{eq:entropie_shift}
\end{eqnarray} 
It remains to use Rothau's Lemma (see Lemma 5.1.4 of \cite{BGL}) that states that for any measure $\mu$ and any constant $a$:
$$\en(g+a)^2\le \en g^2+ 2\integ g^2 \de \mu.$$ 
 Let $p=4$ in Equation \eqref{eq:entropie_shift}. Putting this together with Rothau's Lemma and Theorem \ref{prop:poincare_beta}, we obtain the following:
\begin{align*}
\enb(f^2) & = \integ f^2 \log\left(\dfrac{f^2}{\|f\|^2_{2,\mubb}}\right) \de \mubb  \\
& \leq \enb \left[ (f-<f>_{\beta})^2\right]+2 Var_{\mubb}(f) \vspace{1em}\\
& \leq \eb(f,f) \, \left[ \beta |E| \max_{e \in E} L_e e^{2 \mathcal{D}_G} (1+M\beta) e^{c^\star(U_{\nu}) \beta}  + 4 L^2 (2+\mathcal{D}_G L) e^{2/e}\right] \\
& \leq  C_{\Gamma_G}[1+\beta] e^{c^\star(U_{\nu}) \beta}  \eb(f,f) ,  \\
\end{align*}
where $C_{\Gamma_G}$ is a large enough  constant (independent of $\beta$) that could be made explicit in terms of constants $\mathcal{D}_G$, $|E|$ and $L$ since we trivially have $M \leq \mathcal{D}_G^2$.
\hfill $\square$

\bibliographystyle{ormsv080} 
\bibliography{GGR_01-25.bib} 

\end{document}